\documentclass[11pt,english]{article}
\usepackage{lmodern}
\usepackage{lmodern}
\usepackage[T1]{fontenc}
\usepackage[latin9]{inputenc}
\usepackage{geometry}
\usepackage{xcolor}
\usepackage{babel}
\usepackage{array}
\usepackage{verbatim}
\usepackage{pifont}
\usepackage{bbding}
\usepackage{url}
\usepackage{multirow}
\usepackage{amsmath}
\usepackage{amsthm}
\usepackage{amssymb}
\usepackage[unicode=true,pdfusetitle,
 bookmarks=true,bookmarksnumbered=false,bookmarksopen=false,
 breaklinks=false,pdfborder={0 0 1},backref=false,colorlinks=false]
 {hyperref}
\hypersetup{
 colorlinks,citecolor=blue,filecolor=blue,linkcolor=green,urlcolor=purple}

\makeatletter

\providecommand{\tabularnewline}{\\}

\numberwithin{table}{section}
\theoremstyle{plain}
\newtheorem{thm}{\protect\theoremname}[section]
\theoremstyle{plain}
\newtheorem{prop}[thm]{\protect\propositionname}
\theoremstyle{plain}
\newtheorem{lem}[thm]{\protect\lemmaname}

\usepackage{babel}

\usepackage{blindtext}
\usepackage{url}
\usepackage{enumerate}
\usepackage{setspace}
\usepackage{xcolor}
\usepackage{caption}
\usepackage{multicol}
\usepackage[nocompress]{cite}
\usepackage{wasysym}

\usepackage[ruled,vlined,algosection,linesnumbered]{algorithm2e}
\DontPrintSemicolon
\SetKwProg{Fn}{Function}{:}{}
\SetKwInOut{Require}{Input}
\SetKwInOut{Output}{Output}
\SetKwComment{Comment}{$\triangleright$\ }{}

\SetCommentSty{mycommfont}
\SetArgSty{textnormal}

\SetKwFunction{AIDAL}{AIDAL}
\SetKwFunction{AdapAIDAL}{AdapAIDAL}
\SetKwFunction{ACG}{ACG}

\definecolor{green}{rgb}{0,0.4,0.0} 

\makeatother

\providecommand{\lemmaname}{Lemma}
\providecommand{\propositionname}{Proposition}
\providecommand{\theoremname}{Theorem}

\begin{document}
\global\long\def\tx{\tilde{x}}%
\global\long\def\rn{\mathbb{R}^{n}}%
\global\long\def\R{\mathbb{R}}%
\global\long\def\r{\mathbb{R}}%
\global\long\def\n{\mathbb{N}}%
\global\long\def\c{\mathbb{C}}%
\global\long\def\pt{\mathbb{\partial}}%
\global\long\def\lam{\lambda}%
\global\long\def\argmin{\operatorname*{argmin}}%
\global\long\def\Argmin{\operatorname*{Argmin}}%
\global\long\def\argmax{\operatorname*{argmax}}%
\global\long\def\dom{\operatorname*{dom}}%
\global\long\def\ri{\operatorname*{ri}}%
\global\long\def\diag{\operatorname*{diag}}%
\global\long\def\Diag{\operatorname*{Diag}}%
\global\long\def\inner#1#2{\langle#1,#2\rangle}%
\global\long\def\cConv{\overline{{\rm Conv}}\ }%
\global\long\def\intr{\operatorname*{int}}%
\global\long\def\trc{\operatorname*{tr}}%

\title{An Accelerated Inexact Dampened Augmented Lagrangian Method for Linearly-Constrained
Nonconvex Composite Optimization Problems\thanks{The first author has been supported by (i) the US Department of Energy
(DOE) and UT-Battelle, LLC, under contract DE-AC05-00OR22725, (ii)
the Exascale Computing Project (17-SC-20-SC), a collaborative effort
of the U.S. Department of Energy Office of Science and the National
Nuclear Security Administration, and (iii) the IDEaS-TRIAD Fellowship
(NSF Grant CCF-1740776). The second author was partially supported by ONR Grant N00014-18-1-2077 and AFOSR Grant FA9550-22-1-0088.}}
\date{February 6, 2023 (v1: October 23, 2021; v2: August 12, 2022)}
\author{Weiwei Kong\thanks{Computer Science and Mathematics Division, Oak Ridge National Laboratory,
Oak Ridge, TN, 37830. \protect\href{mailto:wwkong92@gmail.com}{wwkong92@gmail.com}} and Renato D.C. Monteiro\thanks{School of Industrial and Systems Engineering, Georgia Institute of
Technology, Atlanta, GA, 30332-0205. \protect\href{mailto:monteiro@isye.gatech.edu}{monteiro@isye.gatech.edu}}}

\maketitle
 
\begin{abstract}
This paper proposes and analyzes an accelerated inexact dampened augmented
Lagrangian (AIDAL) method for solving linearly-constrained nonconvex
composite optimization problems. Each iteration of the AIDAL method
consists of: (i) inexactly solving a dampened proximal augmented Lagrangian
(AL) subproblem by calling an accelerated composite gradient (ACG)
subroutine; (ii) applying a dampened and under-relaxed Lagrange multiplier
update; and (iii) using a novel test to check whether the penalty
parameter of the AL function should be increased. Under several mild
assumptions involving the dampening factor and the under-relaxation
constant, it is shown that the AIDAL method generates an approximate
stationary point of the constrained problem in ${\cal O}(\varepsilon^{-5/2}\log\varepsilon^{-1})$
iterations of the ACG subroutine, for a given tolerance $\varepsilon>0$.
Numerical experiments are also given to show the computational efficiency
of the proposed method.
\end{abstract}

\section{Introduction}

\label{sec:intro}

This paper presents an accelerated inexact dampened augmented Lagrangian
(AIDAL) method for finding approximate stationary points of the linearly
constrained nonconvex composite optimization (NCO) problem

\begin{equation}
\min_{z}\left\{ \phi(u):=f(z)+h(z):Az=b\right\} ,\label{eq:main_prb}
\end{equation}
where $A$ is a linear operator, $h$ is a proper closed convex and
Lipschitz continuous function with compact domain, and $f$ is a (possibly)
nonconvex differentiable function on the domain of $h$ with a Lipschitz
continuous gradient. More specifically, the AIDAL method is based
on the $\theta$-dampened augmented Lagrangian (AL) function 
\begin{gather}
{\cal L}_{c}^{\theta}(z;p):=\phi(z)+(1-\theta)\left\langle p,Az-b\right\rangle +\frac{c}{2}\|Az-b\|^{2}\quad\forall c>0,\quad\forall\theta\in(0,1),\label{eq:aug_lagr_fn_def}
\end{gather}
and it performs the following updates to generate its $k^{{\rm th}}$
iterate: given $(z_{k-1},p_{k-1})$ and $(\lam,c_{k})$, compute 
\begin{align}
z_{k} & \approx\argmin_{u}\left\{ \lam{\cal L}_{c_{k}}^{\theta}(u;p_{k-1})+\frac{1}{2}\|u-z_{k-1}\|^{2}\right\} ,\label{eq:primal_update}\\
p_{k} & =(1-\theta)p_{k-1}+\chi c_{k}(Az_{k}-b),\label{eq:dual_update}
\end{align}
where $\chi$ is an under-relaxation parameter in $(0,1)$ and $z_{k}$
is a suitably chosen approximate solution of the composite problem
underlying \eqref{eq:primal_update}. In addition, the AIDAL method
introduces a novel approach for updating the penalty parameter $c_{k}$
between iterations and uses an accelerated composite gradient (ACG)
method applied to \eqref{eq:primal_update} obtain the aforementioned
point $z_{k}$.

Under a suitable choice of $\lam$ and the following Slater-like assumption:
\begin{equation}
\exists\bar{z}\in\intr(\dom h)\text{ such that }A\bar{z}=b,\label{eq:strong_Slater}
\end{equation}
where $\intr(\dom h)$ denotes the interior of the domain of $h$,
it is shown that, for any tolerance pair $({\rho},{\eta})\in\r_{++}^{2}$,
the AIDAL method obtains a triple $(\hat{z},\hat{p},\hat{v})$ satisfying
\begin{equation}
\hat{v}\in\nabla f(\hat{z})+\pt h(\hat{z})+A^{*}\hat{p},\quad\|\hat{v}\|\leq{\rho},\quad\|A\hat{z}-b\|\le{\eta}.\label{eq:rho_eta_approx_sol}
\end{equation}
in ${\cal O}(({\eta}^{-5/2}+{\eta}^{-1/2}{\rho}^{-2})\log{\eta}^{-1})$
ACG iterations. Moreover, this iteration complexity is obtained without
requiring that the initial point $z_{0}$ (in the domain of $h$)
be feasible with respect to the linear constraint, i.e., $Az_{0}=b$.
Another contribution from this analysis is that the sequence of Lagrange
multipliers is shown to be bounded by a constant independent
of ${\rho}$ and ${\eta}$.\\

\emph{Related Works.} To condense our discussion, we let $\varepsilon={\rho}={\eta}$
denote a common tolerance parameter and restrict our attention to
works that establish iteration complexity bounds for obtaining approximate
stationary points of \eqref{eq:main_prb}. For an overview of papers
that focus on asymptotic convergence of a proposed method, see the
excellent discussion in \cite[Section 2]{Lin2019}. 

One popular class of methods for obtaining stationary points of \eqref{eq:main_prb}
is the penalty method, which consists of solving a sequence of unconstrained
subproblems containing an objective function that penalizes a violation
of the constraints through a positively weighted penalty term. Papers
\cite{WJRproxmet1,Kong2019} present an ${\cal O}(\varepsilon^{-3})$
iteration complexity of a quadratic penalty method without any regularity
assumptions on the linear constraint. In a follow-up work, paper \cite{WJRVarLam2018}
presents an ${\cal O}(\varepsilon^{-3}\log\varepsilon^{-1})$ iteration
complexity of a similar quadratic penalty method in which its parameters
are chosen in an adaptive and numerically efficient manner.
Paper
\cite{Lin2019} is the first to present a penalty-based method with
an improved complexity of ${\cal O}(\varepsilon^{-5/2}\log\varepsilon^{-1})$
under the assumption that the domain of $h$ is compact and assumption
\eqref{eq:strong_Slater} holds.


Another popular class of methods is the proximal AL (PAL) method,
which primarily consists of the updates in \eqref{eq:primal_update}
and \eqref{eq:dual_update}. The analysis of AL/PAL-based methods
for the case where $\phi$ is convex is already well-established (see,
for example, \cite{LanRen2013PenMet,Aybatpenalty,IterComplConicprog,AybatAugLag,LanMonteiroAugLag,ShiqiaMaAugLag16,zhaosongAugLag18,Patrascu2017,Xu2019}),
so we make no more mention of it here. Instead, we review papers that
present an iteration complexity of an AL/PAL-based method for the
case where $\phi$ is nonconvex. Paper \cite{HongPertAugLag} presents
an ${\cal O}(\varepsilon^{-4})$ iteration complexity\footnote{This method generates prox subproblems of the form $\argmin_{x\in X}\{\lam h(x) + c\|Ax-b\|^2 / 2 + \|x-x_0\|^2 / 2 \}$ and the analysis of \cite{HongPertAugLag}
makes the strong assumption that they can be solved exactly for any $x_0$, $c$, and $\lam$.} of an unaccelerated
PAL method under the strong assumption that the initial point $z_{0}$
is feasible, i.e., $Az_{0}=b$, as well as $\theta\in(0,1]$ and $\chi=1$.
Paper \cite{RJWIPAAL2020} presents ${\cal O}(\varepsilon^{-3}\log\varepsilon^{-1})$
and ${\cal O}(\varepsilon^{-5/2}\log\varepsilon^{-1})$ iteration
complexities of an accelerated inexact PAL method for the general
case and the case where \eqref{eq:strong_Slater} holds, respectively,
and removes the requirement that the initial point be feasible. Papers \cite{Melo2020,kong2020iteration} present an ${\cal O}(\varepsilon^{-3}\log\varepsilon^{-1})$
iteration complexity for the special case of $(\chi,\theta)=(1,0)$,
which corresponds to a full multiplier update under the classical
AL function. Finally, papers \cite{inexactAugLag19}
and \cite{ImprovedShrinkingALM20} respectively establish ${\cal O}(\varepsilon^{-3}\log\varepsilon^{-1})$
and ${\cal O}(\varepsilon^{-5/2}\log\varepsilon^{-1})$ iteration
complexities for nonproximal AL-based methods that perform under-relaxed Lagrange
multiplier updates only when the penalty parameter is updated. 

Aside from penalty and AL/PAL-based methods, we mention few others
that are of interest. Paper \cite{boob2022stochastic} presents an ${\cal O}(\varepsilon^{-3})$
iteration complexity of a primal-dual proximal point scheme for generating
a point \emph{near }an approximate stationary point under some strong
conditions on the initial point. Papers \cite{ErrorBoundJzhang-ZQLuo2020,ADMMJzhang-ZQLuo2020}
present an ${\cal O}(\varepsilon^{-2})$ iteration complexity of a
primal-dual first-order algorithm for solving \eqref{eq:main_prb}
when $h$ is the indicator function of a box (in \cite{ADMMJzhang-ZQLuo2020}),
or more generally, a polyhedron (in \cite{ErrorBoundJzhang-ZQLuo2020}).
Paper \cite{SZhang-Pen-admm} presents an ${\cal O}(\varepsilon^{-6})$
iteration complexity of a penalty-ADMM method that solves an equivalent
reformulation of \eqref{eq:main_prb}, under the assumption that the
initial point $z_{0}$ is feasible, the tolerance $\varepsilon$ is
sufficiently small, and $A$ has full row rank. 
Paper \cite{Li2019} presents an inexact proximal point method applied to the function defined as $\phi(z)$ if $z$ is feasible and $+\infty$ otherwise. 
It can be viewed as an extension to the nonconvex setting
of the proximal point method (PPM) applied to
\eqref{eq:main_prb} and it
 obtains an ${\cal O}(\varepsilon^{-5/2}\log\varepsilon^{-1})$ complexity bound.\\

\emph{Contributions}. We now emphasize how the proposed AIDAL method
improves on other state-of-the-art AL-based works. First, it improves upon the ${\cal O}(\varepsilon^{-3}\log\varepsilon^{-1})$
classic PAL method in \cite{Melo2020} by an ${\cal O}(\varepsilon^{-1/2})$
factor through only a \emph{small} perturbation of the classical multiplier
update and the classical AL function. 
Second, AIDAL chooses its prox stepsize $\lam$ independent of the 
perturbation parameter $\theta$. This is in contrast to the PAL method in \cite{RJWIPAAL2020} which has the undesirable property that its prox stepsize $\lam$ becomes arbitrarily small as $\theta$ approaches zero.
Finally, it differs from the 
nonproximal AL-based method in \cite{ImprovedShrinkingALM20} in two significant
ways: (i) it performs the multiplier update \eqref{eq:dual_update}
after every inexact prox update as opposed to only when the penalty
parameter is updated; and (ii) it chooses a constant under-relaxation
parameter $\chi$ for the update \eqref{eq:dual_update} as opposed
to \cite{ImprovedShrinkingALM20}, which chooses an under-relaxation
parameter that (linearly) tends to zero as the number of penalty parameter
updates increases.  \\


\emph{Organization of the Paper}. Subsection~\ref{subsec:notation}
provides some basic definitions and notation. Section~\ref{sec:main_prb_and_method}
contains two subsections. The first one describes the main problem
of interest and the assumptions made on it, while the second one  presents the AIDAL method
and states its iteration complexity.
Section~\ref{sec:cvg_analysis}
is divided into four subsections. The first one presents some preliminary technical results, the second one presents a bound on an important stationarity residual, the third one proves a bound on the generated Lagrange multipliers, and fourth one one gives the proof
of a key proposition in Section~\ref{sec:main_prb_and_method}. Section~\ref{sec:numerical_experiments}
presents numerical experiments that demonstrate the efficiency of
the AIDAL method. Section~\ref{sec:concl_remarks} gives some concluding
remarks. Finally, the end of the paper contains several important
technical appendices.

\subsection{Basic Notations and Definitions}

\label{subsec:notation}

This subsection presents notation and basic definitions used in this
paper.

Let $\r_{+}$ and $\r_{++}$ denote the set of nonnegative and positive
real numbers, respectively, and let $\r^{n}$ denote the $n$-dimensional
Hilbert space with inner product and associated norm denoted by $\inner{\cdot}{\cdot}$
and $\|\cdot\|$, respectively. The smallest positive singular value
of a nonzero linear operator $Q:\r^{n}\to\r^{l}$ is denoted by $\sigma_{Q}^{+}$.
For a given closed convex set $X\subset\r^{n}$, its boundary is denoted
by $\partial X$ and the distance of a point $x\in\r^{n}$ to $X$
is denoted by ${\rm dist}_{X}(x)$. For any $t>0$, we let $\log_{1}^{+}(t):=\max\{\log t,1\}$ and denote ${\cal O}_{1}={\cal O}(\cdot+1)$. 

The domain of a function $h:\r^{n}\to(-\infty,\infty]$ is the set
$\dom h:=\{x\in\r^{n}:h(x)<+\infty\}$. Moreover, $h$ is said to
be proper if $\dom h\ne\emptyset$. The set of all lower semi-continuous
proper convex functions defined in $\r^{n}$ is denoted by $\cConv\rn$.
The subdifferential of a proper convex function $h:\r^{n}\to(-\infty,\infty]$
is defined by 
\begin{equation}
\partial h(z):=\{u\in\r^{n}:h(z')\geq h(z)+\inner u{z'-z},\quad\forall z'\in\r^{n}\}\label{def:epsSubdiff}
\end{equation}
for every $z\in\r^{n}$. The
normal cone of a closed convex set $C$ at $z\in C$
is defined as 
\[
N_{C}(z):=\{\xi\in\r^{n}:\inner{\xi}{u-z}\leq0,\quad\forall u\in C\}.
\]
If $\psi:\rn \mapsto \r$ is differentiable at $\bar{z}\in\r^{n}$,
then its affine approximation at $\bar{z}$
is given by 
\begin{equation}
\ell_{\psi}(z;\bar{z}):=\psi(\bar{z})+\inner{\nabla\psi(\bar{z})}{z-\bar{z}}\quad\forall z\in\r^{n}.\label{eq:defell}
\end{equation}

\section{Augmented Lagrangian Method}

\label{sec:main_prb_and_method}

This section contains two subsections. The first one precisely describes
the problem of interest and the assumptions underlying it, while the
second one presents the AIDAL method and its corresponding iteration
complexity.

\subsection{Problem of Interest}

\label{subsec:prb_of_interest}

This subsection presents the main problem of interest and the assumptions
underlying it.

Our problem of interest is precisely \eqref{eq:main_prb} where $f$,
$h$, $A$, and $b$ are assumed to satisfy the following assumptions: 
\begin{itemize}
\item[(A1)] $h\in\cConv\rn$ is $K_{h}$-Lipschitz continuous and ${\cal H}:=\dom h$
is compact with diameter $D_{h}:=\sup_{u,z\in{\cal H}}\|u-z\|<\infty$. 
\item[(A2)] $f$ is differentiable function on ${\cal H}$, and there exists
$(m,M)\in\r_{++}^{2}$ satisfying $m\leq M$, such that for every
$u,z\in{\cal H}$, we have 
\begin{gather}
f(u)-\ell_{f}(u;z) \geq -\frac{m}{2}\|u-z\|^{2},\label{eq:curve_cond}\\
\|\nabla f(u)-\nabla f(z)\|\leq M\|u-z\|;\label{eq:Lipschitz_cond}
\end{gather}
\item[(A3)] there exists $\bar{z}\in\intr{\cal H}$ such that $A\bar{z}=b$; 
\item[(A4)] $A\neq0$, ${\cal F}:=\{z\in{\cal H}:Az=b\}\neq\emptyset$, and $\inf_{z\in\rn}\phi(z)>-\infty$. 
\end{itemize}
We now make four remarks about the above assumptions. First, it is
well-known that \eqref{eq:Lipschitz_cond} implies that $|f(u)-\ell_{f}(u;z)|\leq M\|u-z\|^{2}/2$
for every $u,z\in{\cal H}$ and hence that \eqref{eq:curve_cond}
holds with $m=M.$ However, we show that better iteration complexities
can be derived when a scalar $m\ll M$ satisfying \eqref{eq:curve_cond}
is available {\color{purple}(see Theorem~\ref{thm:aidal_compl} and \eqref{eq:spec_acg_compl})}.  Second, \eqref{eq:curve_cond} implies that the function
$f+m\|\cdot\|^{2}/2$ is convex on ${\cal H}$. Third, since ${\cal H}$
is compact by (A1), the image of any continuous $\r^{k}$-valued function,
e.g., $u\mapsto\nabla f(u)$, is bounded. 
Finally, in Appendix~\ref{app:local_necessary}, we show that if $\hat{z}$ is a local minimum of \eqref{eq:main_prb}, then  there exists a multiplier
$\hat{p}$ such that 
\begin{equation}
0\in\nabla f(\hat{z})+\pt h(\hat{z})+A^{*}\hat{p},\quad A\hat{z}=b.\label{eq:stationary_soln}
\end{equation}
In view of the last remark, we say that a triple $(\hat{z},\hat{p},\hat{v})$
is a $({\rho},{\eta})$-stationary point of \eqref{eq:main_prb}
if it satisfies condition \eqref{eq:rho_eta_approx_sol}, which is
clearly a relaxation of \eqref{eq:stationary_soln} for any $({\rho},{\eta})\in\r_{++}^{2}$.

\subsection{AIDAL Method}

\label{subsec:aidal_method}

This section presents the AIDAL method and its corresponding iteration
complexity.

We first state the AIDAL method in Algorithm~\ref{alg:aidal}. Its
main steps are: (i) invoking an ACG algorithm (specifically, Algorithm~\ref{alg:acg}) to implement the update
in \eqref{eq:primal_update}; (ii) computing a ``refined'' pair
$(\hat{p},\hat{v})=(\hat{p}_{k},\hat{v}_{k})$ and point $z$
satisfying the inclusion and (possibly) the inequality in \eqref{eq:rho_eta_approx_sol};
(iii) applying the update in \eqref{eq:dual_update}; and (iv) performing a novel test to determine the next penalty parameter
$c_{k+1}$.

\begin{algorithm}[!htb] 
\caption{Accelerated Inexact Dampened Augmented Lagrangian (AIDAL) Method}
\label{alg:aidal}

\Require{$(m,M)\in\r_{++}^2$ as in (A2), $({\rho},{\eta})\in\r_{++}^{2}$, $(z_{0},p_{0})\in{\cal H}\times A(\r^{n})$, $c_{1} \in \r_{++}$, $\sigma\in(0,1/2]$, and $(\chi,\theta)\in(0,1)^{2}$
satisfying 
\begin{equation}
\color{purple}
(1-\theta)(2-\theta)\chi \le {\theta^{2}}.\label{eq:chi_theta_cond}
\end{equation}
}

\Output{ a triple $(\hat{z},\hat{p},\hat{v})\in {\cal H}\times A(\r^{n})\times\rn$
satisfying \eqref{eq:rho_eta_approx_sol}.}

\BlankLine
\Fn{\AIDAL{$\{m,M\},\{\sigma,\chi,\theta\},\{c_{1},z_{0},p_{0}\},\{{\rho},{\eta}\}$}}{

\texttt{\textcolor{blue}{STEP 0}}\textcolor{blue}{{} (initialization)}

$\lam \gets 1/(2m)$

\For{$k \gets 1,2,...$}{

\texttt{\textcolor{blue}{STEP 1}}\textcolor{blue}{{} (inexact prox update)}:\Comment*[r]{Implement
\eqref{eq:primal_update}}

$L_{k}\gets\lam(M+c_{k}\|A\|^{2})+1$\;

$\psi_{s}^{k}(\cdot)\gets\lam\left[{\cal L}_{c_{k}}^{\theta}(\cdot;p_{k-1})-h(\cdot)\right]+\frac{1}{2}\|\cdot-z_{k-1}\|^{2}$\Comment*[r]{See
\eqref{eq:aug_lagr_fn_def} for the definition of ${\cal L}_{c}^{\theta}(\cdot;\cdot)$}

$({z}_{k},v_{k})\gets$\ACG{$\{\psi_{s}^{k}, \lam h\}, \{L_{k},\frac{1}{2}\},{\sigma},z_{k-1}$}\Comment*[r]{Use
Algorithm~\ref{alg:acg}}

\texttt{\textcolor{blue}{STEP 2}}\textcolor{blue}{{} (termination check)}:\;

$\hat{v}_{k}\gets\frac{1}{\lam}\left[v_{k} + {z}_{k-1} - {z}_{k}\right]$\;

$\hat{p}_{k}\gets(1-\theta)p_{k-1}+c_{k}(A z_{k}-b)$\;

\If{$\|\hat{v}_{k}\|\leq{\rho}$ \textbf{and} $\|A {z}_{k}-b\|\leq {\eta}$}{\Return{$({z}_{k},\hat{p}_{k},\hat{v}_{k})$}\Comment*[r]{Stop
and output}

}

\texttt{\textcolor{blue}{STEP 3}}\textcolor{blue}{{} (multiplier update)}:\Comment*[r]{Implement
\eqref{eq:dual_update}}

$p_{k}\gets(1-\theta)p_{k-1}+\chi c_{k}(Az_{k}-b)$\;

\texttt{\textcolor{blue}{STEP 4}}\textcolor{blue}{{} (penalty parameter
update)}:

$
c_{k+1}\gets\begin{cases}
2c_{k}, & \text{if } \|\hat{v}_k\| \leq {\rho},\\
c_{k}, & \text{otherwise}
\end{cases}
$

} 

} 

\end{algorithm}

Some remarks about Algorithm~\ref{alg:aidal} are in
order. First, its input $z_{0}$ can be any element in ${\cal H}$
and does not necessarily need to be a feasible point, i.e., one satisfying
$Az_{0}=b$. Second, its steps~1 and 3 are respectively the updates
\eqref{eq:primal_update} and \eqref{eq:dual_update}, while its step~4
consists of a test to determine whether the penalty parameter $c_{k}$
should be increased. In particular, the update for \eqref{eq:primal_update}
is obtained by applying the ACG algorithm in Algorithm~\ref{alg:acg}
to the (convex) proximal subproblem
\[
\min_{u\in\rn}\left\{ \lam{\cal L}_{c_{k}}^{\theta}(\cdot;p_{k-1})+\frac{1}{2}\|\cdot-z_{k-1}\|^{2}\right\} 
\]
with an inexactness criterion (see \eqref{eq:approx_acg_soln}) that is a variant of the one considered by the authors in \cite{WJRproxmet1,Kong2019,kong2020iteration, kong2022complexity}.
Third, it performs two kinds of iterations: (i) the ones indexed by
$k$; and (ii) the ones performed by the ACG algorithm every
time it is called in its step~1. To be concise, the former will be
referred to as ``outer'' iterations and the latter as ``inner''
(or ACG) iterations. Finally, it is shown in Lemma~\ref{lem:ext_acg_statn_props}(d) that the triple $(\hat{z}, \hat{p}, \hat{v}) = (z_k,\hat{p}_k,\hat{v}_k)$ satisfies the inclusion in \eqref{eq:rho_eta_approx_sol} for every $k\geq 1$. Hence, if the termination condition in step~3 is satisfied, then AIDAL outputs a $(\rho,\eta)$-stationary point of \eqref{eq:main_prb} (whose definition is given at the end of Subsection~\ref{subsec:prb_of_interest}).

We now present the key properties of the method. To be
concise, we introduce the constants 
\begin{equation}
\begin{gathered}\bar{d}:={\rm dist}_{\pt{\cal H}}(\bar{z}),\quad G_{f}:=\sup_{u\in{\cal H}}\|\nabla f(u)\|,\quad\phi_{*}:=\inf_{u\in\rn}\phi(u),\quad\phi^{*}:=\inf_{u\in{\cal F}}\phi(u),\\
\beta_{\lam}:=\left(\bar{d}+D_{h}\right)\left[K_{h}+G_{f}+\frac{(1+\sigma)D_{h}}{\lam}\right],
\end{gathered}
\label{eq:global_iter_consts}
\end{equation}
where $(D_{h},K_{h},{\cal H})$, $\bar{z}$, and ${\cal F}$ are as
in (A1), (A3), and (A4), respectively. Moreover, we let 
\begin{equation}
{\cal C}_{\ell}:=\left\{ k\in\n: c_{k}=c_{1}2^{\ell-1}\right\} \label{eq:cycle_def}
\end{equation}
denote the $\ell^{{\rm th}}$ \emph{cycle} of AIDAL and, for
simplicity, if the AIDAL terminates at iteration $k$ then
the indices of the last cycle do not extend past $k$.

The first result
presents a bound on the sequence of Lagrange multipliers $\{p_k\}_{k\geq 0}$
computed in step~3 of AIDAL. Its proof, which is given
in Subsection~\ref{subsec:Lagr_mult_bd},
is a generalization of \cite[Proposition~3.12]{Melo2020}, which considers
the case where $(\theta,\chi)=(0,1)$.
\begin{prop}
\label{prop:global_iter_bds}Let $\{p_{i}\}_{i\geq1}$
be generated by the AIDAL method. Then, 
\begin{equation}
\color{purple}
\|p_{k}\|\leq \|p_0\| + \frac{\beta_{\lam}}{\bar{d}\sigma_{A}^{+}} =: B_{p}\quad \forall k\geq 1,
\label{eq:Bp_def}
\end{equation}
where $\bar{d}$ and $\beta_\lam$ are as in \eqref{eq:global_iter_consts}.
\end{prop}

The next result, whose proof is the topic of Subsection~\ref{subsec:alm_props_prf}, describes several properties of AIDAL, including a bound on
the number of inner (or ACG) iterations performed in each outer iteration, a uniform bound on the size of all cycles, and its 
successful termination with the required approximate stationary point of
\eqref{eq:main_prb}.
\begin{prop}
\label{prop:inexact_alm_props} Let $(\lam, c_1, \chi, \theta, \rho, \eta)$ be as in AIDAL, and define the nonnegative scalars
\begin{equation}
\begin{aligned}
B_{\Psi} &:=\phi^{*}-\phi_{*}+\frac{D_{h}^{2}}{\lambda}+
\left(\frac{2 - \theta + 2[2-\theta][1-\theta]}{2\chi^2 c_1}\right)B_p^2, \\
\bar{c}_{{\eta}} &:=\frac{2B_{p}}{\chi{\eta}},
\quad
{\cal T}_{{\rho}}:=\left\lceil 1+\frac{9 B_{\Psi}}{\lambda{\rho}^{2}} \right\rceil.
\end{aligned}
\label{eq:key_props_consts}
\end{equation}
where $B_{p}$, $D_h$, and $(\phi_*, \phi^*)$ are as in Proposition~\ref{prop:global_iter_bds}, assumption (A1), and \eqref{eq:global_iter_consts}, respectively.
Then, the following statements hold about AIDAL: 
\begin{itemize}
\item[(a)] its $k^{\rm th}$ outer iteration performs a number of inner (or ACG) iterations bounded above by
\begin{align}
\color{purple}
\left\lceil 1 + 6\sqrt{L_{k}}\log_{1}^{+} \frac{4L_{k}}{\sigma}\right\rceil, \label{eq:acg_spec_comp}
\end{align}
where $L_k$ is given by step~1 of AIDAL;
\item[(b)] for every $\ell\geq1$, it holds that $|{\cal C}_{\ell}|\leq{\cal T}_{{\rho}}$, and the residual $\hat v_k$ for the last index $k$ of ${\cal C}_\ell$ satisfies $\|\hat{v}_k\|\leq \rho$;
\item[(c)] the last cycle $\bar\ell$ outputs a $({\rho},{\eta})$-stationary point of \eqref{eq:main_prb} and satisfies
$c_k \leq \max\{c_1, 2\bar{c}_\eta\}$ for every $k\in {\cal C}_{\bar{\ell}}$; as a consequence, $\bar{\ell} \leq \max\{1, \log_2(2\bar{c}_\eta/c_1)\}$.
\end{itemize}
\end{prop}

We give some remarks about the above results.
First, Proposition~\ref{prop:global_iter_bds} states that the sequence
of Lagrange multipliers $\{p_{k}\}_{k\geq1}$ generated by the AIDAL
method is bounded by a constant that is independent of the tolerances
${\rho}$ and ${\eta}$. Second, Proposition~\ref{prop:inexact_alm_props}(c)
states that the number of times that the penalty constant $c_{k}$
is doubled during an invocation of the AIDAL method is finite. Finally, Proposition~\ref{prop:inexact_alm_props}(a) shows that the number of the
inner (or ACG) iterations at each outer iteration of AIDAL is independent of the 
tolerances ${\rho}$ and ${\eta}$.

Using Proposition~\ref{prop:inexact_alm_props}, the next result establishes an ${\cal O}({\eta}^{-1/2} {\rho}^{-2} \log {\eta}^{-1})$ total inner (or ACG) iteration complexity for the AIDAL method.
\begin{thm}
\label{thm:aidal_compl} AIDAL stops with a $({\rho},{\eta})$-stationary point of \eqref{eq:main_prb} in a number of inner (or ACG) iterations bounded above by 
\begin{equation}
{\cal O}_{1}\left({\cal T}_{{\rho}}\sqrt{\bar{c}_{{\eta}}  L_{1}}\log_{1}^{+} \frac{\bar{c}_{{\eta}} L_{1}}{\sigma}\right),\label{eq:total_acg_compl}
\end{equation}
where $(\bar{c}_{{\eta}}, {\cal T}_{{\rho}})$ are as \eqref{eq:key_props_consts}, $L_{1}$ is as in step~1 of AIDAL at $k=1$, and $\sigma$ is the inexactness parameter given to the ACG algorithm (Algorithm~\ref{alg:acg}).
\end{thm}

\begin{proof}
For ease of notation, let $(\bar{c}, {\cal T}) = (\bar{c}_{{\eta}}, {\cal T}_{{\rho}})$. In view of Proposition~\ref{prop:inexact_alm_props}(a) and (c),
the total number of inner (or ACG) iterations performed by the method
is on the order of
\begin{equation}
{\cal O}_{1}\left(\sum_{\ell=1}^{\left\lceil \log_{2}\bar{c}\right\rceil }\sum_{j\in{\cal C}_{\ell}}\sqrt{L_{j}}\log_{1}^{+}\frac{L_{j}}{\sigma}\right).\label{eq:cycle_acg_iter}
\end{equation}
To simplify this sum, we first note that if {\color{purple} ${j}\in{\cal C}_{\ell}$},
then the relations $\lam = 1/(2m)$ (from AIDAL) and $m\leq M$
(from assumption (A2)) imply that 
\begin{equation}
L_{j}=\lam\left(M+2^{\ell-1}c_{1}\|A\|^{2}+\lam^{-1}\right)\leq\lam\left(2M+2^{\ell-1}c_{1}\|A\|^{2}\right)\leq2^{\ell}L_{1}.\label{eq:M_c_spec_bd}
\end{equation}
Combining \eqref{eq:M_c_spec_bd} with Proposition~\ref{prop:inexact_alm_props}(b),
it holds that
\begin{align}
\sum_{\ell=1}^{\left\lceil \log_{2}\bar{c}\right\rceil }\sum_{j\in{\cal C}_{\ell}}\sqrt{L_{j}} & \leq{\cal T}\sqrt{L_{1}}\sum_{\ell=1}^{\left\lceil \log_{2}\bar{c}\right\rceil }2^{\ell/2}={\cal T}\sqrt{L_{1}}\cdot\sqrt{2}\left(1+\sqrt{2}\right)\left(2^{\left\lceil \log_{2}\bar{c}\right\rceil /2}-1\right)\nonumber \\
 & \leq4{\cal T}\sqrt{L_{1}}\left(2^{\log_{2}\sqrt{\bar{c}}}\cdot2^{1/2}\right)={\cal O}_{1}\left({\cal {\cal T}}\sqrt{\bar{c}L_{1}}\right).\label{eq:total_acg_aux_bd1}
\end{align}
Moreover, denoting $\bar{\ell}=\left\lceil \log_{2}\bar{c}\right\rceil $,
it follows from \eqref{eq:M_c_spec_bd} that
\begin{align}
\max_{1\leq\ell\leq\left\lceil \log_{2}\bar{c}\right\rceil }\max_{j\in{\cal C}_{\ell}}\left\{ \log_{1}^{+}L_{j}\right\}  & =\log_{1}^{+}\left[\lam\left(M+c_{\bar{\ell}}\|A\|^{2}\right)+1\right]={\cal O}_{1}\left(\log_{1}^{+}\left[\bar{c}L_{1}\right]\right).\label{eq:total_acg_aux_bd2}
\end{align}
The complexity bound in \eqref{eq:total_acg_compl} now follows from
\eqref{eq:total_acg_aux_bd1}, \eqref{eq:total_acg_aux_bd2}, and
\eqref{eq:cycle_acg_iter}. The fact that AIDAL stops with a $(\rho,\eta)$-stationary point of \eqref{eq:main_prb} follows from Proposition~\ref{prop:inexact_alm_props}(c).
\end{proof}
{\color{purple}We now analyze how the complexity bound in \eqref{eq:total_acg_compl} depends on the stepsize $\lam$ and the tolerances $\rho$ and $\eta$. Throughout our discussion, we make the reasonable assumption that the parameter $\chi$ and the initial penalty parameter $c_1$ are not too small in the sense that
$\max\{c_1^{-1}, \chi^{-1}\} = O(1)$. 
In this case, it is easy to see that the quantities $(B_p, B_\Psi, L_1, \bar{c}_\eta,
{\cal T}_\rho)$ in \eqref{eq:Bp_def}, \eqref{eq:key_props_consts}, and step~1 of Algorithm~\ref{alg:aidal} satisfy $B_p = O(1+\lam^{-1})$, $B_\Psi = O([1 + \lam^{-1}]^2)$, $L_1 = O(1 + \lam)$, $\bar{c}_\eta = O([1 +\lam^{-1}]/\eta)$, and ${\cal T}_\rho = O(1 + [1 + \lam^{-1}]^2/[\lam\rho^2])$. Consequently, the bound \eqref{eq:total_acg_compl} is
\begin{equation}
\color{purple}
{\cal O}_{1}\left(\left[1 + \frac{(1+\lambda^{-1})^2}{\lam \rho^{2}}\right]\sqrt{\frac{1+\lambda + \lam^{-1}}{\eta}}\,{\log_{1}^{+}\left[\frac{1+\lam + \lam^{-1}}{{\eta}}\right]}\right). \label{eq:spec_acg_compl}
\end{equation}
Since $\lam^{-1}=O(1)$, the above complexity consists of the sum of two components: $S_1 = O(\lam^{1/2}\eta^{-1/2})$ and $S_2 = O(\lam^{-1/2}\eta^{-1/2}\rho^{-2})$ (ignoring logarithmic terms). 
In general, if the tolerances $\rho$ and $\eta$ are small, then $S_2 \ll S_1$ and choosing larger values of $\lam$ improves the complexity bound in \eqref{eq:spec_acg_compl}. 
Under the assumption that there exists
a constant $m \ll M$
satisfying \eqref{eq:curve_cond},
this observation
justifies the claim made in the paragraph following
assumptions (A1)--(A4),
namely, that AIDAL can benefit
if such $m$ is known;
otherwise, the only option available would be to
be set $\lam$ to the 
much smaller quantity $1/(2M)$.
}

It is also worth mentioning that, the number
of resolvent (or proximal) evaluations of $h$ in AIDAL is on the same
order of magnitude as in \eqref{eq:total_acg_compl} due to the fact
that the ACG algorithm in Appendix~\ref{app:ACG} performs exactly one
resolvent evaluation per ACG iteration. 

\section{Convergence Analysis of the AIDAL Method}

\label{sec:cvg_analysis}

{\color{purple} This section establishes the key properties of the AIDAL method and
contains four subsections. The first one establishes some properties of the ACG call of AIDAL, the second one gives a useful technical bound on the stationarity residuals $\{\hat{v}_i\}$, the third one gives the proof of Proposition~\ref{prop:global_iter_bds}, and the fourth one gives the proof of Proposition~\ref{prop:inexact_alm_props}.}

To avoid repetition, we let 
\[
\{(z_{i},p_{i},v_{i},\hat{p}_{i},\hat{v}_{i}, \psi_{s}^{i},c_{i},L_{i})\}_{i\geq 1},
\]
denote the sequence of iterates generated by the AIDAL method. Moreover, for every
$i\geq1$ and any $(\chi,\theta)\in\r_{++}^{2}$, we make use of the
following useful constants
\begin{equation}
\begin{gathered}a_{\theta}=\theta(1-\theta),\quad b_{\theta}:=(2-\theta)(1-\theta),\quad\alpha_{\chi,\theta}:=\frac{(1-2\chi b_{\theta})-(1-\theta)^{2}}{2\chi},\\
f_{i}:=Az_{i}-b,\quad\Delta p_{i}=p_{i}-p_{i-1},\quad\Delta z_{i}=z_{i}-z_{i-1}.
\end{gathered}
\label{eq:global_alm_consts}
\end{equation}

\subsection{Preliminary Results}

{\color{purple}This subsection establishes two preliminary technical results about  the residuals $v_i$, $\hat{v}_i$, and $f_i$.
It also establishes the iteration-complexity of each ACG call in
step 1 of AIDAL
using the general results derived for this method in Appendix~\ref{app:ACG}.}

\begin{lem}
\label{lem:feas_props}For every $i\geq 1$:
\begin{itemize}
\item[(a)] $f_{i}=\left[p_{i}-(1-\theta)p_{i-1}\right]/(\chi c_{i})$; 
\item[(b)] if $i\geq2$, then $\chi(c_{i}f_{i}-c_{i-1}f_{i-1})=\Delta p_{i}-(1-\theta)\Delta p_{i-1}$; 
\item[(c)] $\|f_{i}\| \leq \color{purple} ({\|p_i\| + (1-\theta)\|p_{i-1}\|})/({\chi c_i})$.
\end{itemize}
\end{lem}

\begin{proof}
(a) This follows from the definition of $f_{i}$ in \eqref{eq:global_alm_consts},
and step~3 of the AIDAL method.

(b) This follows from part (a) and the definition of $\Delta p_{i}$
in \eqref{eq:global_alm_consts}.

(c) Using part (a), the fact that $1-\theta \in [0,1]$, and the triangle inequality, we have
\[
\|f_i\| = \frac{\|p_i - (1-\theta)p_{i-1}\|}{\chi c_i} \leq \frac{\|p_i\| + (1-\theta)\|p_{i-1}\|}{\chi c_i}. \qedhere
\]
\end{proof}

Note that the inequality of {\color{purple}Lemma~\ref{lem:feas_props}(c)} implies the feasibility residual $\|f_i\|$ can be made small by making the penalty parameter sufficiently large {\color{purple} and ensuring that the multipliers $\{p_i\}_{i\geq 1}$ are bounded.}

\begin{lem}
\label{lem:ext_acg_statn_props}
For every $i \geq 1$:
\begin{itemize}
     \item[(a)]  $\psi_s^i(\cdot) - \|\cdot\|_{Q_i}^2/2$ is convex and $\nabla \psi_s^i(\cdot)$ is $L_i$-Lipschitz continuous, where
    $L_i$ is as in step~1 of Algorithm~\ref{alg:aidal} and 
    \begin{equation}
    Q_{i}:= \frac{I}{2}  + c_{i}\lam A^{*}A,\quad\|\cdot\|_{Q_{i}}^2:=\left\langle \cdot , Q_{i}(\cdot)\right\rangle;
    \label{eq:spec_norm}
    \end{equation}
    \item[(b)] {\color{purple}the $i^{\rm th}$ call to Algorithm~\ref{alg:acg} in step~1 of Algorithm~\ref{alg:aidal} stops in a number of ACG iterations bounded above by \eqref{eq:acg_spec_comp}};
    \item[(c)] it holds that
    \[
    v_i \in \pt(\psi_s^i + \lam h)(z_i) = \pt\left(\lam {\cal L}_{c_i}^\theta(\cdot;p_{i-1}) + \frac{1}{2}\|\cdot-z_{i-1}\|^2\right)(z_i), \quad \|v_{i}\|\leq {\sigma}\|\Delta z_i\|;
    \]
    \item[(d)]$\hat{v}_{i}\in\nabla f({z}_{i})+\pt h({z}_{i})+A^{*}\hat{p}_{i}$ and $\|\hat{v}_{i}\| \leq (1+\sigma) \|\Delta z_i\| / \lam$.
\end{itemize}
\end{lem}

\begin{proof} 
{\color{purple} (a) 
First note that
inequality \eqref{eq:curve_cond}
 in Assumption (A2)
 and
 the choice of $\lam = 1/(2m)$ in Algorithm~\ref{alg:aidal} implies that
 $\lam f(\cdot) + \|\cdot-z_{i-1}\|^2/2$
 is $1/2$-strongly convex on ${\cal H}$.
Hence, the convexity assertion follows 
from this observation, the definition of $\psi_s^i$,
 and the definitions of $Q_i$ and $\|\cdot\|_{Q_i}$ in \eqref{eq:spec_norm}. 
On the other hand, the assertion about Lipschitz continuity follows from the definition of $\psi_s^i$ and \eqref{eq:Lipschitz_cond}.}

(b) Using the fact that $L_i\geq 1$ and $\sigma\in(0,1)$, we first observe that for $\mu=1/2$ we have
\begin{align*}
\color{purple}
\frac{4L_i(L_i+\mu)}{\mu\sigma^2} \leq \frac{8L_i(L_i+L_i)^2}{\sigma^2} \leq  \left[\frac{4L_i}{\sigma}\right]^3.
\end{align*}
Then, note that part (a) implies $(\psi_s,\psi_n)=(\psi_s^i,\lam h)$ satisfies assumptions (B1)--(B2) in Appendix~\ref{app:ACG} with $(L,\mu)=(L_i,1/2)$.
The conclusion now follows from step~1 of Algorithm~\ref{alg:aidal}, assumption (A2),
Proposition~\ref{prop:acg_props}(b) with $(L,\mu)=(L_i, 1/2)$, and the above observations.

(c) {\color{purple}Recall that step~1 of AIDAL calls ACG with $(\psi_{s},\psi_{n})=(\psi_{s}^{i},\lam h)$ and $x_0 = z_{i}$. It then follows from  Proposition~\ref{prop:acg_props}(b) that 
\eqref{eq:approx_acg_soln} holds with $(z,v,x_0)=(z_i,v_i,z_{i-1})$ and $(\psi_{s},\psi_{n})=(\psi_{s}^{i},\lam h)$ and $x_0 = z_{i-1}$.
The inclusion and first inequality now follow from the previous observation, the definition of $\psi_s^i$, and the fact that $\psi_s^i + \lambda h$ is convex (see the choice of $\lam$ and assumption (A2)) and, hence, that $\nabla \psi_s^i (\cdot) + \lam \pt h(\cdot) = \pt(\psi_s^i + \lam h)(\cdot)$.}

(d) Using part (c) and
the definitions of $\hat{v}_{i}$, $\hat{p}_i$, and $\psi_s^i$, it holds that 
\begin{align*}
\hat{v}_{i}=\frac{v_{i}+z_{i-1}-z_{i}}{\lambda} & \in \frac{\nabla \psi_s^i(z_i)}{\lam} + \pt h(z_i) +\frac{z_{i-1}-z_{i}}{\lambda}\\
 & =\nabla f(z_{i})+\partial h(z_{i})+(1-\theta)A^{*}p_{i-1}+cA^{*}(Az_{i}-b)\\
 & =\nabla f(z_{i})+\partial h(z_{i})+A^{*}\hat{p}_{i},
\end{align*}
which is the desired inclusion. {\color{purple}For the desired inequality,
 we use part (c), the triangle inequality, and the definition
of $\hat{v}_{i}$ to obtain}
\[
\color{purple}\|\hat{v}_{i}\|=\frac{1}{\lambda}\|v_{i}+z_{i-1}-z_{i}\|\leq\frac{1}{\lambda}\|v_{i}\|+\frac{1}{\lambda}\|\Delta z_{i}\|\leq\frac{1+\sigma}{\lambda}\|\Delta z_{i}\|. \qedhere
\]
\end{proof}

We now make three comments about the above result.
First, statements (a) and (b) of Lemma~\ref{lem:ext_acg_statn_props}
justify
the choice
of $\lam=1/(2m)$.
Second,  $\lam$ could actually have been set to any value $(0,1/m)$
at the expense of more complicated bounds in the resulting analysis.
Third,
in view of the inclusion of {\color{purple}Lemma~\ref{lem:ext_acg_statn_props}(d)} and the definition $f_i$ in \eqref{eq:global_alm_consts},
it follows that
$(z_i,\hat p_i,\hat v_i)$
is a $(\rho,\eta)$-stationary point of \eqref{eq:main_prb} if and only if
$\|\hat v_i\| \le \rho$
and $\|f_i\| \le \eta$.

{\color{purple}In the next subsection, we establish an important bound on the residuals $\{\hat{v}_i\}$ that will be used to show that they tend to zero.}

\subsection{Bounds on the Stationarity Residuals}

\label{subsec:statn_bds}

 This subsection focuses on establishing the following bound on the residuals $\{\hat{v}_i\}_{i\geq 0}$ within cycle
 ${\cal C}_\ell$ for
 any
 $\ell \ge 1$.
 Note that the value of $c_i$ is constant within ${\cal C}_\ell$,
 i.e., there exists $\tilde{c}_\ell>0$ such that
 \begin{equation} \label{eq:tildec}
 c_i = \tilde c_\ell \quad \forall i \in {\cal C}_\ell.
 \end{equation}

\begin{prop}
\label{prop:vi_poten_bd}
For every $\ell\geq 1$ and $j,k\in {\cal C}_\ell$ such that $k\geq j+1$, we have 
\begin{equation}
{\lambda} \sum_{i=j+1}^k \|\hat{v}_{i}\|^{2} \leq 9[\Psi_{j}^{\theta}-\Psi_{k}^{\theta}],
\label{eq:vi_poten_bd}
\end{equation}
where the potential $\Psi_{i}^{\theta}$ is given by
\begin{align}
\Psi_{i}^{\theta} & :={\cal L}_{\tilde{c}_{\ell}}^{\theta}(z_{i};p_{i})-\frac{a_{\theta}}{2\chi \tilde{c}_{\ell}}\|p_{i}\|^{2}+\frac{\alpha_{\chi,\theta}}{4\chi \tilde{c}_{\ell}}\|\Delta p_{i}\|^{2}.\label{eq:Psi_def}
\end{align}
\end{prop}

We start with a technical bound on $\|\hat{v}_i\|$.
\begin{lem}
\label{lem:alm_aux2} For every $i \ge 1$, it holds that
\begin{align}
\begin{aligned}
\frac{\lam}{9} \|\hat v_{i}\|^{2} &\le
\left [
{\cal L}_{c_i}^{\theta}(z_{i-1};p_{i-1}) - {\cal L}_{c_i}^{\theta}(z_{i};p_{i}) +
\frac{a_{\theta}}{2\chi c_i}\left(\|p_{i}\|^{2}-\|p_{i-1}\|^{2}\right) \right] \\
& \qquad + \frac{b_{\theta}}{2\chi c_{i}}\|\Delta p_{i}\|^{2} -
\frac{c_{i}}{2}\|A\Delta z_{i}\|^{2}, 
\end{aligned}
\label{eq:vhat_bd}
\end{align}
where $a_{\theta}$ and $b_{\theta}$ are as in
\eqref{eq:global_alm_consts}.
\end{lem}

\begin{proof}
Let $i \geq 1$ be fixed. 
We first derive a relationship for ${\cal L}_{c_{i}}^{\theta}(z_{i},p_{i})-{\cal L}_{c_{i}}^{\theta}(z_{i},p_{i-1})$. Using the definition of ${\cal L}_{c}^{\theta}$
in \eqref{eq:aug_lagr_fn_def}, the definitions of $\Delta p_{i}$
and $f_{i}$ in \eqref{eq:global_alm_consts}, and Lemma~\ref{lem:feas_props}(a),
we have that 
\begin{align}
{\cal L}_{c_{i}}^{\theta}(z_{i},p_{i})-{\cal L}_{c_{i}}^{\theta}(z_{i},p_{i-1}) & =(1-\theta)\left\langle \Delta p_{i},f_{i}\right\rangle =\left(\frac{1-\theta}{\chi c_{i}}\right)\|\Delta p_{i}\|^{2}+\frac{(1-\theta)\theta}{\chi c_{i}}\left\langle \Delta p_{i},p_{i-1}\right\rangle \nonumber \\
 & =\left(\frac{1-\theta}{\chi c_{i}}\right)\|\Delta p_{i}\|^{2}+\frac{(1-\theta)\theta}{\chi c_{i}}\left(\left\langle p_{i},p_{i-1}\right\rangle -\|p_{i-1}\|^{2}\right)\nonumber \\
 & =\left(\frac{1-\theta}{\chi c_{i}}\right)\|\Delta p_{i}\|^{2}+\frac{(1-\theta)\theta}{\chi c_{i}}\left(-\frac{1}{2}\|\Delta p_{i}\|^{2}+\frac{1}{2}\|p_{i}\|^{2}-\frac{1}{2}\|p_{i-1}\|^{2}\right)\nonumber \\
 & =\frac{b_{\theta}}{2\chi c_{i}}\|\Delta p_{i}\|^{2}+\frac{a_{\theta}}{2\chi c_{i}}\left(\|p_{i}\|^{2}-\|p_{i-1}\|^{2}\right).\label{eq:partial_DeltaL1}
\end{align}
We next derive a bound for ${\cal L}_{c_{i}}^{\theta}(z_{i},p_{i-1})-{\cal L}_{c_{i}}^{\theta}(z_{i-1},p_{i-1})$. 
In view of Lemma~\ref{lem:ext_acg_statn_props}(a) and (c), we first observe that (i) $\lam {\cal L}_{c_{i}}^{\theta}(\cdot,p_{i-1}) + \|\cdot-z_{i-1}\|^2/2 = \psi_s^i (\cdot)+ \lam h(\cdot)$ 
is 1-strongly convex with respect to the $\|\cdot\|_{Q_i}$ norm given in \eqref{eq:spec_norm}, and (ii) $z_i$ is an optimal solution of the function $\psi_s^i (\cdot)+ \lam h(\cdot) - \langle v_i, \cdot \rangle$. 
Combining facts (i)--(ii) above, the definition of $\|\cdot\|_{Q_i}$ in \eqref{eq:spec_norm}, the bound on $\|v_i\|$ in Lemma~\ref{lem:ext_acg_statn_props}(c), the fact that $\sigma\in(0,1/2]$, and the Cauchy-Schwarz inequality, we conclude that
\begin{align}
& {\cal L}_{c_{i}}^{\theta}(z_{i},p_{i-1})-{\cal L}_{c_{i}}^{\theta}(z_{i-1},p_{i-1}) \leq-\frac{1}{2\lam}\|\Delta z_{i}\|_{Q_{i}}^{2}-\frac{1}{2\lam}\|\Delta z_{i}\|^{2}+\frac{1}{\lam}\left\langle v_{i},\Delta z_{i}\right\rangle \nonumber \\
&\le  -\frac{c_{i}}{2}\|A\Delta z_{i}\|^{2} -\frac{3}{4\lam}\|\Delta z_{i}\|^{2} + \frac1\lam \|v_i\| \, \|\Delta z_i\| \nonumber \\
 & \le  -\left(\frac{3 - 4\sigma}{4\lam}\right)\|\Delta z_{i}\|^{2}-\frac{c_{i}}{2}\|A\Delta z_{i}\|^{2}
 \leq 
 {\color{purple}-\frac{1}{4\lam}\|\Delta z_{i}\|^{2}}-\frac{c_{i}}{2}\|A\Delta z_{i}\|^{2} \label{eq:partial_DeltaL2}
\end{align}
The conclusion now follows by summing \eqref{eq:partial_DeltaL1} and \eqref{eq:partial_DeltaL2}, isolating the $\|\Delta z_i\|^2$ term to one side, and using the inequality on $\|\hat{v}_i\|$ in Lemma~\ref{lem:ext_acg_statn_props}(d) with the fact that $(1+\sigma)^2 \leq 9/4$. 
\end{proof}

Note that within a cycle, where the penalty parameters remain constant, the term within the square bracket of the right-hand side of \eqref{eq:vhat_bd} is  telescopic.
Interestingly, the next result shows that the other term on
the right-hand side of \eqref{eq:vhat_bd}
can be telescopically bounded within a fixed
cycle.
It is worth mentioning that the relationship between $\chi$ and $\theta$ in \eqref{eq:chi_theta_cond} plays an important role
in proving this fact.

\begin{lem}
\label{lem:alm_aux3}
For every $i\geq 2$ such that $c_i=c_{i-1}$, it holds that
\begin{equation}
 \frac{b_{\theta}}{2\chi c_{i}}\|\Delta p_{i}\|^{2} - \frac{c_{i}}{2}\|A\Delta z_{i}\|^{2} \leq \frac{\alpha_{\chi,\theta}}{2\chi c_{i}}\left(\|\Delta p_{i-1}\|^{2}-\|\Delta p_{i}\|^{2}\right), \label{eq:tech_vhat_sub_bd}
\end{equation}
where $b_\theta$ and $\alpha_{\chi,\theta}$ are as in \eqref{eq:global_alm_consts}.
\end{lem}
\begin{proof}
Let $i\geq2$ be an index where $c_{i}=c_{i-1}$ and observe that \eqref{eq:chi_theta_cond}
implies $2\chi b_{\theta}\leq\theta^{2}$. Moreover, define 
\[
\widehat{\Delta}p_{i}:=\Delta p_{i}-(1-\theta)\Delta p_{i-1}
\]
and observe that Lemma~\ref{lem:tech_ADeltaZ_DeltaP} with $(\tau,a,b)=(\chi b_{\theta},\Delta p_{i},\Delta p_{i-1})$
implies that 
\begin{equation}
\frac{1}{\chi}\|\widehat{\Delta}p_{i}\|^{2}\geq2b_{\theta}\|\Delta p_{i}\|^{2}+\alpha_{\chi,\theta}\left(\|\Delta p_{i}\|^{2}-\|\Delta p_{i-1}\|^{2}\right).\label{eq:hat_DeltaP_aux1}
\end{equation}
Using Lemma~\ref{lem:feas_props}(b), the fact that $c_{i}=c_{i-1}$,
and \eqref{eq:hat_DeltaP_aux1}, we then have 
\begin{align*}
\frac{c_{i}}{2}\|A\Delta z_{i}\|^{2} 
& =\frac{\|\chi c_{i}A\Delta z_{i}\|^{2}}{2\chi^{2}c_{i}}
=\frac{\|\chi(c_{i}f_{i}-c_{i-1}f_{i-1})\|^{2}}{2\chi^{2}c_{i}} 
=\frac{1}{2\chi c_{i}}\left[\frac{1}{\chi}\|\widehat{\Delta}p_{i}\|^{2}\right] \\
& \geq\frac{1}{2\chi c_{i}}\left[b_{\theta}\|\Delta p_{i}\|^{2}+\alpha_{\chi,\theta}\left(\|\Delta p_{i}\|^{2}-\|\Delta p_{i-1}\|^{2}\right)\right],
\end{align*}
from which \eqref{eq:tech_vhat_sub_bd} immediately follows.
\end{proof}

Combining \eqref{eq:vhat_bd} and \eqref{eq:tech_vhat_sub_bd}, it is easy to see that the sum of the residuals $\{\|\hat{v}_i\|^2\}_{i\geq 1}$ residuals is bounded above by a telescopic sum when the indices are in a cycle. Let us now use this fact to prove Proposition~\ref{prop:vi_poten_bd}.

\begin{proof}[Proof of Proposition~\ref{prop:vi_poten_bd}] Let $\ell\geq 1$ and
$j,k\in{\cal C}_\ell$ be given and assume that
$i \in \{j+1,\ldots,k\}$.
Then, it follows from
\eqref{eq:tildec} that
$c_{i-1}=c_i=\tilde c_\ell$.
This observation together
Lemmas~\ref{lem:alm_aux2} and \ref{lem:alm_aux3}
then imply that
\begin{align*}
\frac{\lambda}{9}\|\hat{v}_{i}\|^{2} 
&  \leq{\cal L}_{c_{i-1}}^{\theta}(z_{i-1};p_{i-1})-{\cal L}_{c_{i}}^{\theta}(z_{i};p_{i})+\frac{a_{\theta}}{2\chi c_{i}}\left(\|p_{i}\|^{2}-\|p_{i-1}\|^{2}\right)+\frac{\alpha_{\chi,\theta}}{4\chi c_{i}}\left(\|\Delta p_{i-1}\|^{2}-\|\Delta p_{i}\|^{2}\right)\\
& = {\cal L}_{\tilde{c}_\ell}^{\theta}(z_{i-1};p_{i-1})-{\cal L}_{\tilde{c}_\ell}^{\theta}(z_{i};p_{i})+\frac{a_{\theta}}{2\chi \tilde{c}_\ell}\left(\|p_{i}\|^{2}-\|p_{i-1}\|^{2}\right)+\frac{\alpha_{\chi,\theta}}{4\chi \tilde{c}_\ell}\left(\|\Delta p_{i-1}\|^{2}-\|\Delta p_{i}\|^{2}\right)\\
 & =\Psi_{i-1}^{\theta}-\Psi_{i}^{\theta},
\end{align*}
where the second identity is due to
the definition of $\Psi_i^\theta$ in \eqref{eq:Psi_def}.
The conclusion now follows by summing the above inequality from $i=j+1$ to $k$.
\end{proof}

{\color{purple}
One of the goals of the following two subsections is to show that the potential $\Psi_i^\theta$ in \eqref{eq:Psi_def} can be bounded
by a constant
that does not depend on $c_i$.
A key  step in this direction
is given by
Proposition~\ref{prop:global_iter_bds}
which states
that
the Lagrange multiplier $p_i$
can also be bounded by a constant that does not depend on $c_i$.
The goal of the next subsection
is to prove this proposition.}


\subsection{Proof of Proposition~\ref{prop:global_iter_bds}}

\label{subsec:Lagr_mult_bd}


We start by  presenting two well-known technical results. The proof of the first one can be found, for example, in \cite[Lemma 1.2]{Goncalves2019}. 
\begin{lem}
\label{lem:lower_spectral_bd}For every $S\in\r^{m\times n}$ and $u\in{\rm Im}\ S$,
we have $\sigma_{S}^{+}\|u\|\leq\|Su\|$. 
\end{lem}

The proof of the next result can be found in 
\cite[Lemma~3.10]{Melo2020}.

\begin{lem}
\label{lem:topo_slater}Suppose $\psi\in\cConv\rn$ is $K_{\psi}$-Lipschitz
continuous with finite diameter $D_{\psi}$. Then, for every
$y,\bar{y}\in\dom h$ and $\xi\in\pt\psi(y)$, we have
\[
\|\xi\|{\rm dist}_{\pt(\dom\psi)}(\bar{y})\leq\left[{\rm dist}_{\pt(\dom\psi)}(\bar{y})+\|y-\bar{y}\|\right]K_{\psi}+\left\langle \xi,y-\bar{y}\right\rangle.
\]
\end{lem}

The next two results closely follow the ones in \cite[Section 3]{Melo2020}.

\begin{lem}
\label{lem:iter_props}Define the scalars
\begin{equation}
\xi_{k}:=\hat{v}_{k}-\nabla f(z_{k})-A^{*}\hat{p}_{k}\quad\forall k\geq 1.\label{eq:xi_def}
\end{equation}
Then, the following statements hold for every $k\geq1$:
\begin{itemize}
\item[(a)] $\xi_{k}\in\pt h(z_{k})$;
\item[(b)] it holds that 
\[
\|\hat{p}_{k}\|\leq\frac{1}{\sigma_{A}^{+}}\left[\|\xi_{k}\|+ G_f +\frac{(1+\sigma)D_h}{\lam}\right],
\]
\color{purple} where $G_f$ and $D_h$ are as in \eqref{eq:global_iter_consts} and assumption (A1), respectively.
\end{itemize}
\end{lem}

\begin{proof}
(a) This follows immediately from 
{\color{purple}Lemma~\ref{lem:ext_acg_statn_props}(d)}
and the definition of $\xi_{i}$. 

(b) Using the definitions of $\xi_{i}$ and $G_f$, the triangle inequality, {\color{purple}part
(a)}, and Lemma~\ref{lem:lower_spectral_bd} with $S=A^{*}$ and $u=\hat{p}_{k}$
yields 
\begin{align*}
\|\hat{p}_{k}\| & \leq\frac{\|A^{*}\hat{p}_{k}\|}{\sigma_{A}^{+}}=\frac{\|\hat{v}_{k}-\nabla f(z_{k})-\xi_{k}\|}{\sigma_{A}^{+}} \leq \frac{\|\xi_{k}\|+\|\nabla f(z_{k})\|+ \|\hat{v}_k\|}{\sigma_{A}^{+}} \\
& \leq\frac{1}{\sigma_{A}^{+}}\left[\|\xi_{k}\|+ \|\nabla f(z_{k})\| +\frac{(1+\sigma)\|\Delta z_k\|}{\lam}\right] \leq \frac{1}{\sigma_{A}^{+}}\left[\|\xi_{k}\|+ G_f +\frac{(1+\sigma)D_h}{\lam}\right]. \qedhere
\end{align*}
\end{proof}

\begin{lem}
\label{lem:tech_dual_ineqs}Let $(\beta_{\lam},\bar{d})$ be as in
\eqref{eq:global_iter_consts}. Then, the following statements hold
for every $(\chi,\theta)\in(0,1)^{2}$ and $k\geq1$: 
\begin{itemize}
\item[(a)] $\|p_{k}\|\leq\chi\|\hat{p}_{k}\|+(1-\chi)(1-\theta)\|p_{k-1}\|;$
\item[(b)] $c_{k}^{-1}\|\hat{p}_{k}\|^{2}+\bar{d}\sigma_{A}^{+}\|\hat{p}_{k}\|\leq c_{k}^{-1}\left(1-\theta\right)\left\langle \hat{p}_{k},p_{k-1}\right\rangle +\beta_{\lam}.$
\end{itemize}
\end{lem}

\begin{proof}
(a) Using the definitions of $p_{k}$ and $\hat{p}_{k}$ with the
triangle inequality yields 
\[
\|p_{k}\|=\|\chi\hat{p}_{k}+(1-\chi)(1-\theta)p_{k-1}\|\leq\chi\|\hat{p}_{k}\|+(1-\chi)(1-\theta)\|p_{k-1}\|.
\]

(b) Let $\xi_{k}$, $(G_{f},\bar{d})$, and $D_{h}$ be as in \eqref{eq:xi_def},
\eqref{eq:global_iter_consts}, and assumption (A1), respectively.
Using Lemma~\ref{lem:iter_props}(a), the definition of $\bar{d}$,
and Lemma~\ref{lem:topo_slater} with $(\psi,K_{\psi},D_{\psi})=(h,K_{h},D_{h})$
and $(y,\bar{y},\varepsilon)=({z}_{k},\bar{z},\delta_{k})$, we
have that 
\begin{align}
\bar{d}\|\xi_{k}\| & \leq(\bar{d}+D_{h})K_{h}+\left\langle \xi_{k},z_{k}-\bar{z}\right\rangle.\label{eq:tech_dual_ineq1}
\end{align}
Moreover, the definitions of $\hat{p}_{k}$ and $\xi_{k}$, the
fact that $z_{k},\bar{z}\in{\cal H}$ and $A\bar{z}=b$, and the Cauchy-Schwarz
inequality imply that
\begin{align}
\left\langle \xi_{k},z_{k}-\bar{z}\right\rangle  & 
=\left\langle \hat{v}_{k}-\nabla f(z_{k})-A^{*}\hat{p}_{k},z_{k}-\bar{z}\right\rangle \nonumber \\
 & \leq\left(\|\hat{v}_k\|+\|\nabla f(z_{k})\|\right)\|z_{k}-\bar{z}\|-\left\langle \hat{p}_{k},Az_{k}-b\right\rangle \nonumber \\
 & \leq\left[\frac{(1+\sigma)D_{h}}{\lam}+G_{f}\right]D_{h}+\left(\frac{1-\theta}{c_{k}}\right)\left\langle \hat{p}_{k},p_{k-1}\right\rangle -\frac{1}{c_{k}}\|\hat{p}_{k}\|^{2}.\label{eq:tech_dual_ineq2}
\end{align}
Using Lemma~\ref{lem:iter_props}(b), \eqref{eq:tech_dual_ineq1},
\eqref{eq:tech_dual_ineq2}, and the definition of $\beta_{\lam}$ in \eqref{eq:global_iter_consts}, we thus conclude
that
\begin{align*}
\frac{1}{c_{k}}\|\hat{p}_{k}\|^{2}+\bar{d}\sigma_{A}^{+}\|\hat{p}_{k}\| 
&\leq\frac{1}{c_{k}}\|\hat{p}_{k}\|^{2}+\bar{d}\|\xi_{k}\|+\left[G_{f}+\frac{(1+\sigma)D_{h}}{\lam}\right]\bar{d}\\
 & \color{purple} 
 \leq 
  \frac{1}{c_{k}}\|\hat{p}_{k}\|^{2} + (\bar{d} + D_h)K_h + \langle\xi_k, z_k-\bar{z} \rangle
  +\left[G_{f}+\frac{(1+\sigma)D_{h}}{\lam}\right]\bar{d}
 \\
 & \color{purple}
 \leq 
 \left(\frac{1-\theta}{c_{k}}\right)\left\langle \hat{p}_{k},p_{k-1}\right\rangle +
 \left[K_h + G_{f}+\frac{(1+\sigma)D_{h}}{\lam}\right](\bar{d} + D_h)
 \\
 &
  \color{purple}
 = \left(\frac{1-\theta}{c_{k}}\right)\left\langle \hat{p}_{k},p_{k-1}\right\rangle +\beta_{\lam}.
 \qedhere
\end{align*}
\end{proof}

We are now ready to give the proof of Proposition~\ref{prop:global_iter_bds}.

\begin{proof}[Proof of Proposition~\ref{prop:global_iter_bds}]
We proceed by induction on $k$. Since $B_{p}\geq\|p_{0}\|$, the
desired bound trivially holds for $k=0$. Assume now that $\|p_{k}\|\leq B_{p}$
holds for some $k\geq0$. If $\|\hat{p}_{k+1}\|=0$, then clearly
\[
\|p_{k+1}\|\leq\chi\|\hat{p}_{k+1}\|+(1-\chi)(1-\theta)\|p_{k}\|=(1-\chi)(1-\theta)B_{p}\leq B_{p},
\]
so suppose that $\|\hat{p}_{k+1}\|>0$. Using Lemma~\ref{lem:tech_dual_ineqs}(b),
the Cauchy-Schwarz inequality, and the induction hypothesis we have
that
\begin{align*}
& \left[\bar{d}+\frac{1}{c_{k+1}\sigma_{A}^{+}}\|\hat{p}_{k+1}\|\right]\|\hat{p}_{k+1}\| \leq\frac{1}{\sigma_{A}^{+}}\left[\left(\frac{1-\theta}{c_{k+1}}\right)\left\langle \hat{p}_{k+1},p_{k}\right\rangle +\beta_{\lam}\right]\\
 & \leq \frac{\beta_{\lam}}{\sigma_{A}^+} + \frac{(1-\theta)\|p_{k}\|\cdot\|\hat{p}_{k+1}\|}{c_{k+1}\sigma_{A}^{+}} 
 {\color{purple} \leq 
 \frac{\beta_{\lam}}{\sigma_{A}^+} + \frac{\|\hat{p}_{k+1}\|B_p}{c_{k+1}\sigma_{A}^{+}}
\leq\left[\bar{d}+\frac{1}{c_{k+1}\sigma_{A}^{+}}\|\hat{p}_{k+1}\|\right]B_{p},}
\end{align*}
and, hence, that $\|\hat{p}_{k+1}\|\leq B_{p}$. Combining this bound
with the induction hypothesis, we finally conclude that
\[
\|p_{k+1}\|\leq\chi\|\hat{p}_{k+1}\|+(1-\chi)(1-\theta)\|p_{k}\|\leq B_{p}. \qedhere
\]
\end{proof}

\subsection{Proof of Proposition~\ref{prop:inexact_alm_props}}

\label{subsec:alm_props_prf}

{\color{purple}
Recall that Proposition~\ref{prop:vi_poten_bd} in Subsection~\ref{subsec:statn_bds} gives a bound on $\sum_{i=j+1}^k\|\hat{v}_i\|^2$ in \eqref{eq:vi_poten_bd}. 
The first part of this subsection further refines
\eqref{eq:vi_poten_bd} to show that its right-hand side is bounded by a constant that does not depend on the constant $\tilde{c}_\ell$ in \eqref{eq:tildec}. The following result provides a key step in this direction.
}




\begin{lem}
\label{lem:L_c_bds}
For every $i\geq1$, it holds that
\begin{equation}
\phi_{*}-\left(\frac{1-\theta}{2 \chi c_1}\right) B_{p}^{2} 
\leq \Psi_i^\theta
\leq \phi^{*}+\frac{D_{h}^{2}}{\lam} + \left(\frac{1 + 2b_{\theta}}{2\chi^2 c_1}\right)B_p^2,\label{eq:Lagr_bds}
\end{equation}
where  $(\phi_{*},\phi^{*})$, $b_\theta$, and $D_{h}$ are as
in \eqref{eq:global_iter_consts}, \eqref{eq:global_alm_consts}, and assumption (A1), respectively.
\end{lem}

\begin{proof}
Let $i\geq1$. Using Proposition~\ref{prop:global_iter_bds}, the definitions of ${\cal L}_{c}^{\theta}(\cdot,\cdot)$, $\Psi_j^\theta$,
$\phi_{*}$, and $B_{p}$, and the fact that $\chi \in (0,1)$, we have 
\begin{align*}
\Psi_i^\theta & \geq  {\cal L}_{c_{i}}^{\theta}(z_{i};p_{i}) - \frac{a_\theta}{2\chi c_i}\|p_i\|^2 =\phi(z_{i}) + (1-\theta)\left\langle p_{i},Az_{i}-b\right\rangle +\frac{c_{i}}{2}\|Az_{i}-b\|^{2} - \frac{a_\theta}{2\chi c_i}\|p_i\|^2 \\
 & \geq\phi_{*}+\frac{1}{2}\left\Vert \left(\frac{1-\theta}{\sqrt{c_{i}}}\right)p_{i}+\sqrt{c_{i}}(Az_{i}-b)\right\Vert ^{2}-\frac{(1-\theta)^{2}}{2c_{i}}\|p_{i}\|^{2} 
 - \frac{a_\theta}{2\chi c_i}\|p_i\|^2 \\ 
 & \geq\phi_{*}- \left[\frac{(1-\theta)^{2} + a_\theta}{2\chi c_{i}}\right] B_{p}^{2} \geq \phi_{*}- \left(\frac{1-\theta}{2 \chi c_1}\right) B_{p}^{2},
\end{align*}
which is the desired lower bound in \eqref{eq:Lagr_bds}. For the
upper bound, let an arbitrary $u\in{\cal F}$ be given.
Using the fact that $Au=b$ and
$u \in {\cal H}$, the definitions of ${\cal L}_c^\theta(\cdot,\cdot)$ and $D_h$ ,  Lemma~\ref{lem:ext_acg_statn_props}(c), and the Cauchy-Schwarz inequality,
we conclude that
\begin{align*}
\lambda{\cal L}_{c_{i}}^{\theta}(z_{i};p_{i-1}) 
& \overset{\text{Lemma}~\ref{lem:ext_acg_statn_props}(c)}{\leq}
\lambda{\cal L}_{c_{i}}^{\theta}(u;p_{i-1})+\frac{1}{2}\|u-z_{i-1}\|^{2}-\frac{1}{2}\|\Delta z_{i}\|^{2}-\left\langle v_{i},u-z_{i}\right\rangle \\
& \overset{u\in{\cal F}}{\leq}
\lambda\phi(u)+\frac{1}{2} D_h^2 + {\color{purple}\| v_{i}\|} D_h
 \overset{\text{Lemma}~\ref{lem:ext_acg_statn_props}(c)}{\leq}
\lambda\phi(u)+ \left( \frac{1}{2} + \sigma \right) D_h^2.
\end{align*}
Taking the infimum of the above bound over $u\in {\cal F}$ and using the fact that $\sigma \in (0, 1/2]$, we 
thus have ${\cal L}_{c_{i}}^{\theta}(z_{i};p_{i-1})\leq\phi^{*}+D_{h}^{2}/\lam$.
This inequality, \eqref{eq:partial_DeltaL1}, the fact that $\chi\in(0,1)$, Proposition~\ref{prop:global_iter_bds},
 and the relation $(a+b)^{2}\leq2a^{2}+2b^{2}$
for every $a,b\in\r$, then imply that
\begin{align*}
\Psi_{i}^{\theta} 
& = 
{\cal L}_{c_{i}}^{\theta}(z_{i};p_{i}) -\frac{a_\theta}{2\chi c_i}\|p_i\|^2 + \frac{\alpha_{\chi,\theta}}{4 \chi c_i}\|\Delta p_i\|^2 
\leq 
{\cal L}_{c_{i}}^{\theta}(z_{i};p_{i-1})+\left(\frac{2b_{\theta} + \alpha_{\chi,\theta}}{4\chi c_{i}}\right)\|\Delta p_{i}\|^{2}\\
 & \leq
 \phi^{*}+\frac{D_{h}^{2}}{\lam}+\left(\frac{2b_{\theta} + \alpha_{\chi,\theta}}{2\chi c_{i}}\right)(\|p_{i}\|^2 + \|p_{i-1}\|^2) 
 \leq 
 \phi^{*}+\frac{D_{h}^{2}}{\lam} + \left(\frac{1 + 2b_{\theta}}{2\chi^2 c_1}\right)B_p^2,
\end{align*}
which is the desired upper bound in \eqref{eq:Lagr_bds}. 
\end{proof}

The result below follows as a consequence of Proposition~\ref{prop:vi_poten_bd} and Lemma~\ref{lem:L_c_bds}.
\begin{lem}
\label{lem:main_resid_bd} For every $\ell\geq1$ and $j,k \in {\cal C}_{\ell}$ such that $j<k$,
there exists $i\in\{j+1,...,k\}$ satisfying
\begin{equation}
\lambda\|\hat{v}_{i}\|^{2}\leq\frac{9 B_{\Psi}}{k-j},\label{eq:main_resid_bd}
\end{equation}
where $B_{\Psi}$ is as in \eqref{eq:key_props_consts}.
\end{lem}

\begin{proof}
Using the first bound of \eqref{eq:Lagr_bds} with $i=k$ and the second bound of \eqref{eq:Lagr_bds} with $i=j$, we first have that 
\[
\Psi_{j}^\theta-\Psi_{k}^\theta 
\leq 
\phi^{*}-\phi_{*}+\frac{D_{h}^{2}}{\lambda}+
\left(\frac{2 - \theta + 2b_{\theta}}{2\chi^2 c_1}\right)B_p^2 
=
B_\Psi,
\]
where $B_\Psi$ is as in \eqref{eq:key_props_consts}.
Using the above bound and Proposition~\ref{prop:vi_poten_bd}, it follows that
\begin{align*}
  \lam (k-j)\min_{j+1\leq i\leq k}\|\hat v_{i}\|^{2} 
  & \leq
  \lam \sum_{i=j+1}^{k}\|\hat v_{i}\|^{2} 
 \leq
  9 \left(\Psi_{j}^{\theta}-\Psi_{k}^{\theta}\right) 
 \leq 9 B_{\Psi},
\end{align*}
which implies the existence of some $i\in\{j+1,...,k\}$ satisfying
the bound on $\|\hat{v}_{i}\|$ in \eqref{eq:main_resid_bd}.
\end{proof}

We are now ready to give the proof of Proposition~\ref{prop:inexact_alm_props}. 

\begin{proof}[Proof of Proposition~\ref{prop:inexact_alm_props}.]

{\color{purple}(a) This follows immediately from Lemma~\ref{lem:ext_acg_statn_props}(b).}

(b) The fact that the last index $k$ of a cycle ${\cal C}_\ell$ satisfies $\|\hat{v}_k\|\leq \rho$ follows immediately from steps 2--3 of AIDAL.
Now, let $\ell\geq1$ be fixed and define $j:=\inf\{i : i\in{\cal C}_{\ell}\}$
and $k:=j+{\cal T}_\rho-1$. If $k \not\in {\cal C}_\ell$ then $|{\cal C}_{\ell}|\leq k-j+1={\cal T}_\rho$. On the other hand, if $k\in {\cal C}_\ell$ then Lemma~\ref{lem:main_resid_bd} and the definition of ${\cal T}_\rho$ in \eqref{eq:key_props_consts} imply that there exists $i\in\{j+1,...,k\}$
such that
\[
\|\hat{v}_{i}\|^2 \leq \frac{9 B_{\Psi}}{\lam ({\cal T}_\rho-1)} \leq {\rho}^2.
\]
Since every cycle stops when $\|\hat{v}_i\|\leq{\rho}$, we conclude that $i=k=\sup\{i : i \in {\cal C}_\ell\}$ and, hence, $|{\cal C}_\ell| = k - j + 1 = {\cal T}_\rho$.

(c) Let $\bar{c}={\bar{c}}_\eta$. We first establish the bound on $c_k$. If AIDAL stops in the first cycle, then the bound on $c_k$ follows immediately. Assume now that there is more than one cycle and suppose, for the sake of contradiction, that there exists a cycle $\ell \geq 2$ such that $c_k > 2 \bar{c}$ for every $k \in {\cal C}_\ell$, and let $k'$ denote the last index in ${\cal C}_{\ell -1}$. In view steps~3 of AIDAL, we then have $c_{k'} > \bar{c}$. Using the previous bound, the definition of $\bar{c}=\bar{c}_\eta$ in \eqref{eq:key_props_consts}, Lemma~\ref{lem:feas_props}(c), and Proposition~\ref{prop:global_iter_bds}, we also have
\[
\|A x_{k'} - b\| = \|f_{k'}\| \leq {\color{purple}\frac{\|p_{k'}\| + (1-\theta)\|p_{k'-1}\|}{\chi c_{k'}}} \leq \frac{2B_p}{\chi c_{k'}} \leq 
\frac{2B_p}{\chi \bar{c}} \leq \eta.
\]
However, since $\|\hat{v}_{k'}\|\leq \rho$ from part (b), this is impossible because termination would have occurred at the end of cycle $\ell-1$. Hence, $c_k \leq \max\{c_1, 2\bar{c}\}$. Since $c_k = 2^{\bar{\ell}-1} c_1$ for every $k\in {\cal C}_{\bar \ell}$, the bound on $\bar{\ell}$ is immediate. 
Moreover, it follows from parts (a)--(b) and the fact that $\bar{\ell}$ is finite  that AIDAL always stops in step~2. Hence, using the termination condition in step~2 and the inclusion in Lemma~\ref{lem:ext_acg_statn_props}(d), we conclude that the output of AIDAL is a $(\rho,\eta)$-stationary point of \eqref{eq:main_prb}.
\end{proof}

\section{Numerical Experiments}

\label{sec:numerical_experiments}

This section examines the performance of the AIDAL method for solving
problems of the form given in \eqref{eq:main_prb}. It contains four
subsections. The first three contain the following problem classes:
(i) a class of linearly-constrained quadratic programming problems
considered in \cite{WJRproxmet1}; (ii) the sparse principal component
analysis (PCA) problem in \cite{NIPS2014_5615}; and (iii) a class
of linearly-constrained quadratic matrix problems considered in \cite{WJRVarLam2018,kong2020iteration}. The last subsection gives a 
few comments about the results.

Before proceeding with the results, we describe the implementation details of our algorithms and the setup of our experiments. These include specific parameter choices, special modifications, and added heuristics.

We first discuss the three implementation
of the AIDAL method, labeled rADL0, rADL1, and tADL1
considered in this section. Broadly speaking, tADL1 is an implementation of the theoretical version of AIDAL in Algorithm~\ref{alg:aidal}, while rADL0 and rADL1 are implementations of an adaptive/relaxed version of AIDAL in Algorithm~\ref{alg:adap_aidal}. In particular, the adaptive version of AIDAL introduces a novel line search scheme for adaptively choosing the prox parameter $\lam$ in AIDAL (for further details, see the discussion in Appendix~\ref{app:adapt_aidal}).
In terms of parameters, each AIDAL implementation uses 
$p_{0}=0$, $c_{1}=\max\{1,M/\|A\|^{2}\}$, and $\sigma=0.3$ for every outer iteration of the method. However, rADL0 chooses $(\chi,\theta, \lam_0)=(1,0,10)$ with a heuristic choice of $\alpha_{\chi,\theta}=0$ and $a_\theta=1$ in the definition of $\Psi_i^\theta$, while rADL1 and tADL1 choose $(\chi,\theta)=(1/6, 1/2)$ and $\lam_0=10$ for rADL.  Note that rADL0
uses parameters that do not satisfy \eqref{eq:chi_theta_cond},
but work well in practice. 

Besides the above AIDAL implementations, we also use four other
methods as benchmarks. The first one, named iALM, is an implementation
of the inexact proximal augmented Lagrangian method of \cite{ImprovedShrinkingALM20}
in which: (i) its key parameters are
\[
\sigma=5,\quad\beta_{0}=\max\left\{ 1,\frac{\max\{m,M\}}{\|A\|^{2}}\right\} ,\quad w_{0}=1,\quad\boldsymbol{y}^{0}=0,\quad\gamma_{k}=\frac{\left(\log2\right)\|Ax^{1}\|}{(k+1)\left[\log(k+2)\right]^{2}},
\]
for every $k\geq1$; and (ii) the starting point given to the $k^{{\rm th}}$
APG call is set to be $\boldsymbol{x}^{k-1}$, which is the prox center
for the $k^{{\rm th}}$ prox subproblem. The second one, named IPL,
is an implementation of the inexact proximal augmented Lagrangian
method of \cite[Section 5]{kong2020iteration} where: (i) $c_{k}$
is doubled in its step~4 rather than quintupled; and (ii) $\sigma=0.3$.
The third one, named QP, is a practical modification of the quadratic
penalty method of \cite{WJRproxmet1} in which: (i) each ACG subproblem
in step~1 of the AIPP method is stopped when the condition
\[
\|u_{j}\|+2\eta_{j}\leq\sigma\|x_{0}-x_{j}+u_{j}\|^{2}
\]
holds; and (ii) it uses the parameters $\sigma=0.3$ and $c=\max\{1,M/\|A\|^{2}\}$. The fourth and last one, named RQP, is an instance of the relaxed
quadratic penalty method of \cite{WJRVarLam2018} in which: (i) it
uses the AIPPv1 variant described in \cite[Section 6]{WJRVarLam2018}
with the parameters $(\theta,\tau)=(4,10[\lam_0 M+1])$ and $\lam_0=10$;
and (ii) it uses the initial penalty parameter $c_{1}=\max\{1,M/\|A\|^{2}\}$.
It is also worth mentioning that every method except the iALM replaces
its ACG prox subproblem solver by a more practical FISTA variant whose
key iterates are as described in \cite{MontSvaiter_fista} and whose
main stepsize parameter is adaptively estimated by a line search subroutine
described in \cite[Algorithm 5.2.1]{kong2021thesis}.

We now give some comments about the benchmark algorithms.
First, iALM differs from the other tested methods
in that it uses an ACG variant with a termination criterion that is
different from the one in \eqref{eq:approx_acg_soln} and/or its relaxation.
Second, the main difference
between the AIDAL variants and IAIPAL methods is in how they decide when to
double $c_{k}$, i.e., step~4 of Algorithm~\ref{alg:aidal}. In
particular, the condition used in the IAIPAL method depends on both $\sigma$
and $k$ whereas the condition in the AIDAL variants do not. Finally, 
QP-AIPP is the only method that
can be run without requiring any regularity conditions on the linear
constraint and without assuming that $D_{h}<\infty$. 
 In Table~\ref{tab:method_comp}, we summarize the adaptivity of the above methods in terms of the adaptivity of the 
curvature constants $M$ and $m$ in assumption (A2). In particular, we consider the adaptivity
of $m$ to be equivalent to the adaptivity of the prox stepsize $\lam$.

\begin{table}[tbh]
\begin{centering}
\begin{tabular}{>{\centering}p{2.2cm}|>{\centering}m{0.8cm}>{\centering}m{0.8cm}>{\centering}m{0.8cm}>{\centering}m{0.8cm}>{\centering}m{0.8cm}>{\centering}m{0.8cm}>{\centering}m{0.8cm}}
{\footnotesize{}Properties} & {\scriptsize{}rADL0} & {\scriptsize{}rADL1} & {\scriptsize{}tADL1} & {\scriptsize{}iALM} & {\scriptsize{}IPL} & {\scriptsize{}QP} & {\scriptsize{}RQP}\tabularnewline
\hline 
{\scriptsize{}Estimates $M$} & {\footnotesize{}\Checkmark{}} & {\footnotesize{}\Checkmark{}} & {\footnotesize{}\Checkmark{}} & {\footnotesize{}\ding{56}} & {\footnotesize{}\Checkmark{}} & {\footnotesize{}\Checkmark{}} & {\footnotesize{}\Checkmark{}}\tabularnewline
{\scriptsize{}Estimates $m$} & {\footnotesize{}\Checkmark{}} & {\footnotesize{}\Checkmark{}} & {\footnotesize{}\ding{56}} & {\footnotesize{}\ding{56}} & {\footnotesize{}\ding{56}} & {\footnotesize{}\ding{56}} & {\footnotesize{}\Checkmark{}}\tabularnewline
\end{tabular}
\par\end{centering}
\caption{The
first (resp. second) row indicates whether a line search is used to
estimate the curvature constant $M$ (resp. $m$) in assumption (A2)
for a prox subproblem. Note that estimation of $m$ is equivalent
to estimation of the prox stepsize $\protect\lam$. \label{tab:method_comp}}
\end{table}

For a linear operator $A$, a proper lower semicontinuous convex function
$h$, a function $f$ satisfying assumptions (A2)--(A4), a tolerance
pair $({\rho},{\eta})\in\r_{++}^{2}$, and an initial point
$z_{0}\in\dom h$, each of the methods of this section seeks a pair
$([\hat{z},\hat{p}],\hat{v})$ satisfying
\begin{gather}
\begin{gathered}\hat{v}\in\nabla f(\hat{z})+\pt h(\hat{z})+A^{*}\hat{p},\quad\frac{\|\hat{v}\|}{\|\nabla f(z_{0})\|+1}\leq{\rho},\quad\frac{\|A\hat{z}-b\|}{\|Az_{0}-b\|+1}\le{\eta}.\end{gathered}
\label{eq:term_lin_constr}
\end{gather}
In particular, the quadratic programming and matrix problem experiments consider $({\rho},{\eta})=(10^{-3},10^{-3})$, while the sparse PCA experiments consider $({\rho},{\eta})=(10^{-4},10^{-4})$.  Moreover, defining $c_0$ to be the initial penalty parameter and $n_i$ to be the number of outer iterations with $c=c_0 2^{i}$, we also report the following metrics:
\[
c_{\rm wavg} := \frac{\sum_{i\geq 0} n_i \cdot c_0 2^{i}}{\sum_{i \geq 0} n_i} , \quad c_{\max} := \text{ final penalty parameter } c.
\]
All experiments are implemented in MATLAB 2020b and are run on
Linux 64-bit machines, each containing Xeon E5520 processors and at
least 8 GB of memory. Furthermore, the bold numbers in each of the
tables of this section indicate the method that performed the most
efficiently for a given benchmark, e.g., runtime or (innermost) iteration
count. Finally, it is worth mentioning that the code for replicating
these experiments is freely available online\footnote{See \url{https://github.com/wwkong/nc_opt/tree/master/tests/papers/aidal}.}.

\subsection{Linearly-Constrained Quadratic Programming}

\label{subsec:qvp}

Given a pair of dimensions $(l,n)\in\n^{2}$, scalar pair $(\alpha_{1},\alpha_{2})\in\r_{++}^{2}$,
matrices $A,B,C\in\r^{l\times n}$, positive diagonal matrix $D\in\r^{n\times n}$,
and vector pair $(b,d)\in\r^{l}\times\r^{l}$, this subsection considers
the following linearly-constrained quadratic programming (LCQP) problem:
\begin{align*}
\min_{z}\  & \frac{\alpha_{1}}{2}\|Cz-d\|^{2}-\frac{\alpha_{2}}{2}\|DBz\|^{2}\\
\text{s.t.}\  & Az=b,\quad z\in\Delta_{n},
\end{align*}
where $\Delta_{n}=\{z\in\r_{+}^{n}:\sum_{i=1}^{n}z_{i}=1\}$ denotes
the $n$-dimensional simplex.

We now describe the experiment parameters for the instances considered.
First, the dimensions are set to $(l,n)=(10,50)$ and \emph{all} of
the entries in $A$, $B$, and $C$ are nonzero. Second, the entries
of $A,B,C,b$, and $d$ (resp., $D$) are generated by sampling from
the uniform distribution ${\cal U}[0,1]$ (resp., ${\cal U}[1,1000]$).
Third, the initial starting point $z_{0}$ is generated by sampling
a random vector $\tilde{z}_{0}$ from ${\cal U}^{2}[0,1]$ and setting
$z_{0}=\tilde{z}_{0}/\|\tilde{z}_{0}\|$. Fourth, using the well-known
fact that $\|z\|\leq1$ for every $z\in\Delta_{n}$, the auxiliary
parameters for the iALM are $B_{i}=\|a_{i}\|$, $L_{i}=0$, and $\rho_{i}=0$, for every $i$,  
where $a_{i}$ is the $i^{{\rm th}}$ row of $A$. Finally, the composite
form of the problem is 
\[
f(z)=\frac{\alpha_{1}}{2}\|Cz-d\|^{2}-\frac{\alpha_{2}}{2}\|DBz\|^{2},\quad h(z)=\delta_{\Delta_{n}}(z),
\]
and each problem instance uses a scalar pair $(\alpha_{1},\alpha_{2})\in\r_{++}^{2}$
so that $M=\lambda_{\max}(\nabla^{2}f)$
is a particular value given in the table below and $m=-M/3$.

We now present the numerical results for this set of problem instances
in Table~\ref{tab:qvp_iter_run} and Table~\ref{tab:qvp_penalty}.

\begin{table}[tbh]
\begin{centering}
\makebox[\textwidth][c]{%
\begin{tabular}{>{\centering}m{0.7cm}|>{\centering}p{0.5cm}>{\centering}p{0.5cm}>{\centering}p{0.5cm}>{\centering}p{0.5cm}>{\centering}p{0.5cm}>{\centering}p{0.5cm}>{\centering}p{0.5cm}|>{\centering}p{0.4cm}>{\centering}p{0.4cm}>{\centering}p{0.4cm}>{\centering}p{0.4cm}>{\centering}p{0.4cm}>{\centering}p{0.4cm}>{\centering}p{0.4cm}}
{\scriptsize{}$M$} & \multicolumn{7}{c|}{{\scriptsize{}Iteration Count}} & \multicolumn{7}{c}{{\scriptsize{}Runtime (seconds)}}\tabularnewline
\cline{2-15} \cline{3-15} \cline{4-15} \cline{5-15} \cline{6-15} \cline{7-15} \cline{8-15} \cline{9-15} \cline{10-15} \cline{11-15} \cline{12-15} \cline{13-15} \cline{14-15} \cline{15-15} 
 & {\tiny{}rADL0} & {\tiny{}rADL1} & {\tiny{}tADL1} & {\tiny{}iALM} & {\tiny{}IPL} & {\tiny{}QP} & {\tiny{}RQP} & {\tiny{}rADL0} & {\tiny{}rADL1} & {\tiny{}tADL1} & {\tiny{}iALM} & {\tiny{}IPL} & {\tiny{}QP} & {\tiny{}RQP}\tabularnewline
\hline 
{\scriptsize{}$10^{2}$} & \textbf{\tiny{}958} & {\tiny{}1196} & {\tiny{}6910} & {\tiny{}11498} & {\tiny{}26256} & {\tiny{}20473} & {\tiny{}2455} & \textbf{\tiny{}2.0} & {\tiny{}2.5} & {\tiny{}14.0} & {\tiny{}13.8} & {\tiny{}53.4} & {\tiny{}37.9} & {\tiny{}4.6}\tabularnewline
{\scriptsize{}$10^{3}$} & {\tiny{}2538} & {\tiny{}2807} & {\tiny{}7307} & {\tiny{}12669} & {\tiny{}25846} & {\tiny{}20354} & \textbf{\tiny{}2261} & {\tiny{}5.2} & {\tiny{}5.7} & {\tiny{}15.8} & {\tiny{}17.1} & {\tiny{}53.9} & {\tiny{}38.2} & \textbf{\tiny{}4.2}\tabularnewline
{\scriptsize{}$10^{4}$} & \textbf{\tiny{}856} & {\tiny{}2624} & {\tiny{}7307} & {\tiny{}12729} & {\tiny{}25846} & {\tiny{}20497} & {\tiny{}2710} & \textbf{\tiny{}1.7} & {\tiny{}5.4} & {\tiny{}15.2} & {\tiny{}15.8} & {\tiny{}53.0} & {\tiny{}38.4} & {\tiny{}5.0}\tabularnewline
{\scriptsize{}$10^{5}$} & \textbf{\tiny{}908} & {\tiny{}2649} & {\tiny{}7322} & {\tiny{}12743} & {\tiny{}25846} & {\tiny{}20311} & {\tiny{}4571} & \textbf{\tiny{}1.8} & {\tiny{}5.3} & {\tiny{}14.7} & {\tiny{}15.0} & {\tiny{}52.6} & {\tiny{}38.5} & {\tiny{}8.8}\tabularnewline
\textbf{\scriptsize{}$10^{6}$} & \textbf{\tiny{}1045} & {\tiny{}2514} & {\tiny{}7322} & {\tiny{}12744} & {\tiny{}25846} & {\tiny{}20313} & {\tiny{}7889} & \textbf{\tiny{}2.1} & {\tiny{}5.2} & {\tiny{}15.2} & {\tiny{}15.8} & {\tiny{}60.0} & {\tiny{}39.9} & {\tiny{}14.8}\tabularnewline
\end{tabular}}
\par\end{centering}
\caption{Innermost iteration counts and runtimes for LCQP problems.\label{tab:qvp_iter_run}}
\end{table}

\begin{table}[tbh]
\begin{centering}
\makebox[\textwidth][c]{%
\begin{tabular}{>{\centering}m{0.7cm}|>{\centering}m{0.5cm}>{\centering}m{0.5cm}>{\centering}m{0.5cm}>{\centering}m{0.5cm}>{\centering}m{0.5cm}>{\centering}m{0.5cm}>{\centering}m{0.5cm}|>{\centering}m{0.4cm}>{\centering}m{0.4cm}>{\centering}m{0.4cm}>{\centering}m{0.4cm}>{\centering}m{0.4cm}>{\centering}m{0.4cm}>{\centering}m{0.4cm}}
{\scriptsize{}$M$} & \multicolumn{7}{c|}{{\scriptsize{}$c_{\max}$}} & \multicolumn{7}{c}{{\scriptsize{}$c_{{\rm wavg}}/c_{\max}$}}\tabularnewline
\cline{2-15} \cline{3-15} \cline{4-15} \cline{5-15} \cline{6-15} \cline{7-15} \cline{8-15} \cline{9-15} \cline{10-15} \cline{11-15} \cline{12-15} \cline{13-15} \cline{14-15} \cline{15-15} 
 & {\tiny{}rADL0} & {\tiny{}rADL1} & {\tiny{}tADL1} & {\tiny{}iALM} & {\tiny{}IPL} & {\tiny{}QP} & {\tiny{}RQP} & {\tiny{}rADL0} & {\tiny{}rADL1} & {\tiny{}tADL1} & {\tiny{}iALM} & {\tiny{}IPL} & {\tiny{}QP} & {\tiny{}RQP}\tabularnewline
\hline 
{\scriptsize{}$10^{2}$} & \textbf{\tiny{}6E+1} & {\tiny{}2E+3} & {\tiny{}2E+3} & {\tiny{}3E+3} & {\tiny{}3E+5} & {\tiny{}4E+3} & {\tiny{}4E+3} & {\tiny{}0.10} & {\tiny{}0.15} & {\tiny{}0.02} & {\tiny{}0.02} & {\tiny{}0.75} & {\tiny{}0.20} & {\tiny{}0.08}\tabularnewline
{\scriptsize{}$10^{3}$} & \textbf{\tiny{}2E+3} & {\tiny{}4E+4} & {\tiny{}4E+4} & {\tiny{}3E+4} & {\tiny{}3E+6} & {\tiny{}4E+4} & {\tiny{}4E+4} & {\tiny{}0.12} & {\tiny{}0.14} & {\tiny{}0.01} & {\tiny{}0.02} & {\tiny{}0.75} & {\tiny{}0.19} & {\tiny{}0.10}\tabularnewline
{\scriptsize{}$10^{4}$} & \textbf{\tiny{}2E+4} & {\tiny{}4E+5} & {\tiny{}4E+5} & {\tiny{}3E+5} & {\tiny{}3E+7} & {\tiny{}4E+5} & {\tiny{}4E+5} & {\tiny{}0.18} & {\tiny{}0.13} & {\tiny{}0.01} & {\tiny{}0.02} & {\tiny{}0.75} & {\tiny{}0.20} & {\tiny{}0.13}\tabularnewline
{\scriptsize{}$10^{5}$} & \textbf{\tiny{}2E+5} & {\tiny{}4E+6} & {\tiny{}4E+6} & {\tiny{}3E+6} & {\tiny{}3E+8} & {\tiny{}4E+6} & {\tiny{}4E+6} & {\tiny{}0.18} & {\tiny{}0.13} & {\tiny{}0.01} & {\tiny{}0.02} & {\tiny{}0.75} & {\tiny{}0.19} & {\tiny{}0.14}\tabularnewline
\textbf{\scriptsize{}$10^{6}$} & \textbf{\tiny{}2E+6} & {\tiny{}4E+7} & {\tiny{}4E+7} & {\tiny{}3E+7} & {\tiny{}3E+9} & {\tiny{}4E+7} & {\tiny{}4E+7} & {\tiny{}0.18} & {\tiny{}0.13} & {\tiny{}0.01} & {\tiny{}0.02} & {\tiny{}0.75} & {\tiny{}0.19} & {\tiny{}0.15}\tabularnewline
\end{tabular}}
\par\end{centering}
\caption{Penalty parameter statistics for LCQP problems.\label{tab:qvp_penalty}}
\end{table}

It is worth mentioning that we also attempted to add the sProxALM method of \cite{ADMMJzhang-ZQLuo2020,ErrorBoundJzhang-ZQLuo2020} to our list of benchmark methods with its penalty parameter set to $\Gamma=10$ and all other parameters set as in \cite[Algorithm~2]{ErrorBoundJzhang-ZQLuo2020}.
However, for every problem instance, sProxALM failed to obtain a solution as in \eqref{eq:term_lin_constr} under a generous time limit of 3600 seconds, so we have excluded its addition to the results above. Note that we did not test sProxALM on the other numerical experiments because their settings did not fall into settings considered by \cite{ADMMJzhang-ZQLuo2020,ErrorBoundJzhang-ZQLuo2020} (i.e., where the composite function $h$ needs to be the indicator function for a polyhedral set). 
Also, contrary to our AIDAL implementations, \cite{ADMMJzhang-ZQLuo2020,ErrorBoundJzhang-ZQLuo2020} does not provide a concrete way of choosing the parameters (adaptively or otherwise) of sProxALM to ensure its convergence. 

\subsection{Sparse PCA}
\label{subsec:spca}
Given integer $k$, positive scalar pair $(\nu,b)\in\r_{++}^{2}$,
and matrix $\Sigma\in S_{+}^{n}$, this subsection considers the following
sparse principal component analysis (SPCA) problem:
\begin{align*}
\min_{\Pi,\Phi}\  & \left\langle \Sigma,\Pi\right\rangle _{F}+\sum_{i,j=1}^{n}q_{\nu}(\Phi_{ij})+\nu\sum_{i,j=1}^{n}|\Phi_{ij}|\\
\text{s.t.}\  & \Pi-\Phi=0,\quad(\Pi,\Phi)\in{\cal F}^{k}\times\r^{n\times n},
\end{align*}
where ${\cal F}^{k}=\{z\in S_{+}^{n}:0\preceq z\preceq I,\trc M=k\}$
denotes the $k$--Fantope and $q_{\nu}(\cdot)+\nu|\cdot|$ is the
minimax concave penalty (MCP) function given by
\[
q_{\nu}(t):=\begin{cases}
-t^{2}/(2b), & \text{if }|t|\leq b\nu,\\
b\nu^{2}/2-\nu|t|, & \text{if }|t|>b\nu,
\end{cases}\quad\forall t\in\r.
\]
Note that the effective domain of this problem is unbounded, and hence,
only the QP method is guaranteed to converge to an approximate stationary
point in general. 

We now describe the experiment parameters for the instances considered.
First, the scalar parameters are chosen to be $(\nu,b)=(100,0.005)$.
Second, the matrix $\Sigma$ is generated according to an eigenvalue
decomposition $\Sigma=P\Lambda P^{T}$, based on a parameter pair
$(s,k)$, where $k$ is as in the problem description and $s$ is
a positive integer. In particular, we choose $\Lambda=(100,1,...,1)$,
the first column of $P$ to be a sparse vector whose first $s$ entries
are $1/\sqrt{s}$, and the other entries of $P$ to be sampled randomly
from the standard Gaussian distribution. Third, the initial starting
point is $(\Pi_{0},\Phi_{0})=(D_{k},0)$ where $D_{k}$ is a diagonal
matrix whose first $k$ entries are 1 and whose remaining entries
are 0. Fourth, the curvature parameters for each problem instance
are $m=M=1/b$ and $k$ is fixed at $k=1$. Fifth, for the iALM, we
make the following parameter choices based on a relaxed (but unverified)
assumption that its generated iterates lie in ${\cal F}_{k}\times{\cal F}_{k}$:
$B_{i}=1$, $L_{i}=0$, and $\rho_{i}=0$ for all $i$.
Sixth, the composite form of the problem is 
\begin{gather*}
f(\Pi,\Phi)=\left\langle \Sigma,\Pi\right\rangle _{F}+\sum_{i,j=1}^{n}q_{\nu}(\Phi_{ij}),\quad h(\Pi,\Phi)=\delta_{{\cal F}^{k}}(\Pi)+\nu\sum_{i,j=1}^{n}|\Phi_{ij}|,\\
A(\Pi,\Phi)=\Pi-\Phi,\quad b=0,
\end{gather*}
and each problem instance considers a different value of $s$.

We now present the numerical results for this set of problem instances
in Tables~\ref{tab:spca_iter_run} and \ref{tab:spca_penalty}. 

\begin{table}[tbh]
\begin{centering}
\makebox[\textwidth][c]{%
\begin{tabular}{>{\centering}m{0.7cm}|>{\centering}p{0.5cm}>{\centering}p{0.5cm}>{\centering}p{0.5cm}>{\centering}p{0.5cm}>{\centering}p{0.5cm}|>{\centering}p{0.5cm}>{\centering}p{0.5cm}>{\centering}p{0.5cm}>{\centering}p{0.5cm}>{\centering}p{0.5cm}}
{\scriptsize{}$s$} & \multicolumn{5}{c|}{{\scriptsize{}Iteration Count}} & \multicolumn{5}{c}{{\scriptsize{}Runtime (seconds)}}\tabularnewline
\cline{2-11} \cline{3-11} \cline{4-11} \cline{5-11} \cline{6-11} \cline{7-11} \cline{8-11} \cline{9-11} \cline{10-11} \cline{11-11} 
 & {\tiny{}rADL0} & {\tiny{}iALM} & {\tiny{}IPL} & {\tiny{}QP} & {\tiny{}RQP} & {\tiny{}rADL0} & {\tiny{}iALM} & {\tiny{}IPL} & {\tiny{}QP} & {\tiny{}RQP}\tabularnewline
\hline 
{\scriptsize{}$5$} & \textbf{\tiny{}394} & {\tiny{}44952} & {\tiny{}2779} & {\tiny{}22559} & {\tiny{}2990} & \textbf{\tiny{}3.0} & {\tiny{}139.2} & {\tiny{}17.0} & {\tiny{}118.1} & {\tiny{}16.6}\tabularnewline
{\scriptsize{}$10$} & \textbf{\tiny{}403} & {\tiny{}47373} & {\tiny{}2646} & {\tiny{}19984} & {\tiny{}2983} & \textbf{\tiny{}2.7} & {\tiny{}143.1} & {\tiny{}14.8} & {\tiny{}103.8} & {\tiny{}15.8}\tabularnewline
{\scriptsize{}$15$} & \textbf{\tiny{}398} & {\tiny{}45552} & {\tiny{}2628} & {\tiny{}20126} & {\tiny{}2996} & \textbf{\tiny{}2.4} & {\tiny{}138.2} & {\tiny{}15.1} & {\tiny{}103.8} & {\tiny{}16.6}\tabularnewline
\end{tabular}}
\par\end{centering}
\caption{Innermost iteration counts and runtimes for SPCA problems.\label{tab:spca_iter_run}}
\end{table}

\begin{table}[tbh]
\begin{centering}
\makebox[\textwidth][c]{%
\begin{tabular}{>{\centering}m{0.7cm}|>{\centering}p{0.5cm}>{\centering}p{0.5cm}>{\centering}p{0.5cm}>{\centering}p{0.5cm}>{\centering}p{0.5cm}|>{\centering}p{0.5cm}>{\centering}p{0.5cm}>{\centering}p{0.5cm}>{\centering}p{0.5cm}>{\centering}p{0.5cm}}
{\scriptsize{}$s$} & \multicolumn{5}{c|}{{\scriptsize{}$c_{\max}$}} & \multicolumn{5}{c}{{\scriptsize{}$c_{{\rm wavg}}/c_{\max}$}}\tabularnewline
\cline{2-11} \cline{3-11} \cline{4-11} \cline{5-11} \cline{6-11} \cline{7-11} \cline{8-11} \cline{9-11} \cline{10-11} \cline{11-11} 
 & {\tiny{}rADL0} & {\tiny{}iALM} & {\tiny{}IPL} & {\tiny{}QP} & {\tiny{}RQP} & {\tiny{}rADL0} & {\tiny{}iALM} & {\tiny{}IPL} & {\tiny{}QP} & {\tiny{}RQP}\tabularnewline
\hline 
{\scriptsize{}$5$} & \textbf{\tiny{}6E+3} & {\tiny{}4E+6} & {\tiny{}3E+5} & {\tiny{}4E+6} & {\tiny{}2E+6} & {\tiny{}0.57} & {\tiny{}0.03} & {\tiny{}0.33} & {\tiny{}0.04} & {\tiny{}0.09}\tabularnewline
{\scriptsize{}$10$} & \textbf{\tiny{}6E+3} & {\tiny{}4E+6} & {\tiny{}3E+5} & {\tiny{}4E+6} & {\tiny{}2E+6} & {\tiny{}0.57} & {\tiny{}0.03} & {\tiny{}0.28} & {\tiny{}0.03} & {\tiny{}0.09}\tabularnewline
{\scriptsize{}$15$} & \textbf{\tiny{}6E+3} & {\tiny{}4E+6} & {\tiny{}3E+5} & {\tiny{}4E+6} & {\tiny{}2E+6} & {\tiny{}0.57} & {\tiny{}0.03} & {\tiny{}0.35} & {\tiny{}0.03} & {\tiny{}0.09}\tabularnewline
\end{tabular}}
\par\end{centering}
\caption{Penalty parameter statistics for SPCA problems.\label{tab:spca_penalty}}
\end{table}

\subsection{Linearly-Constrained Quadratic Matrix Problem}

\label{subsec:qmp}

Given a pair of dimensions $(l,n)\in\n^{2}$, scalar pair $(\alpha_{1},\alpha_{2})\in\r_{++}^{2}$,
linear operators ${\cal A}:S_{+}^{n}\mapsto\r^{l}$ , ${\cal B}:S_{+}^{n}\mapsto\r^{n}$,
and ${\cal C}:S_{+}^{n}\mapsto\r^{l}$ defined by
\[
\left[{\cal A}(z)\right]_{i}=\left\langle A_{i},z\right\rangle ,\quad\left[{\cal B}(z)\right]_{j}=\left\langle B_{j},z\right\rangle ,\quad\left[{\cal C}(z)\right]_{i}=\left\langle C_{i},z\right\rangle ,
\]
for matrices $\{A_{i}\}_{i=1}^{l},\{B_{j}\}_{j=1}^{n},\{C_{i}\}_{i=1}^{l}\subseteq\r^{n\times n}$,
positive diagonal matrix $D\in\r^{n\times n}$, and vector pair $(b,d)\in\r^{l}\times\r^{l}$,
this subsection considers the following linearly-constrained quadratic
matrix (LCQM) problem:
\begin{align*}
\min_{z}\  & \frac{\alpha_{1}}{2}\|{\cal C}(z)-d\|^{2}-\frac{\alpha_{2}}{2}\|D{\cal B}(z)\|^{2}\\
\text{s.t.}\  & {\cal A}(z)=b,\quad z\in P_{n},
\end{align*}
where $P_{n}=\{z\in S_{+}^{n}:\trc z=1\}$ denotes the $n$-dimensional
spectraplex.

We now describe the experiment parameters for the instances considered.
First, the dimensions are set to $(l,n)=(20,100)$ and only 1.0\%
of the entries of the submatrices $A_{i},B_{j},$ and $C_{i}$ are
nonzero. Second, the entries of $A_{i},B_{j},C_{i},b$, and $d$ (resp.,
$D$) are generated by sampling from the uniform distribution ${\cal U}[0,1]$
(resp., ${\cal U}[1,1000]$). Third, the initial starting point $z_{0}$
is a random point in $S_{+}^{n}$. More specifically, three unit vectors
$\nu_{1},\nu_{2},\nu_{3}\in\r^{n}$ and three scalars $e_{1},e_{2},e_{2}\in\r_{+}$
are first generated by sampling vectors $\tilde{\nu}_{i}\sim{\cal U}^{n}[0,1]$
and scalars $\tilde{d}_{i}\sim{\cal U}[0,1]$ and setting $\nu_{i}=\tilde{\nu}_{i}/\|\tilde{\nu}_{i}\|$
and $e_{i}=\tilde{e}_{i}/(\sum_{j=1}^{3}\tilde{e}_{i})$ for $i=1,2,3$.
The initial iterate for the first subproblem is then set to $z_{0}=\sum_{i=1}^{3}e_{i}\nu_{i}\nu_{i}^{T}$.
Fourth, using the well-known fact that $\|z\|_{F}\leq1$ for every
$z\in P_{n}$, the auxiliary parameters for the iALM are
\[
B_{i}=\|A_{i}\|_{F},\quad L_{i}=0,\quad\rho_{i}=0\quad\forall i\geq1.
\]
Finally, the composite form of the problem is 
\[
f(z)=\frac{\alpha_{1}}{2}\|{\cal C}(z)-d\|^{2}-\frac{\alpha_{2}}{2}\|D{\cal B}(z)\|^{2},\quad h(z)=\delta_{P_{n}}(z),\quad A(z)={\cal A}(z),
\]
and each problem instance uses a scalar pair $(\alpha_{1},\alpha_{2})\in\r_{++}^{2}$
so that $M=\lambda_{\max}(\nabla^{2}f)$
is a particular value given in the table below and $m=-M/4$.

We now present the numerical results for this set of problem instances
in Tables~\ref{tab:qmp_iter_run} and \ref{tab:qmp_penalty}. 

\begin{table}[tbh]
\begin{centering}
\makebox[\textwidth][c]{%
\begin{tabular}{>{\centering}m{0.7cm}|>{\centering}p{0.6cm}>{\centering}p{0.6cm}>{\centering}p{0.6cm}>{\centering}p{0.6cm}>{\centering}p{0.6cm}|>{\centering}p{0.6cm}>{\centering}p{0.6cm}>{\centering}p{0.6cm}>{\centering}p{0.6cm}>{\centering}p{0.6cm}}
{\scriptsize{}$M$} & \multicolumn{5}{c|}{{\scriptsize{}Iteration Count}} & \multicolumn{5}{c}{{\scriptsize{}Runtime (seconds)}}\tabularnewline
\cline{2-11} \cline{3-11} \cline{4-11} \cline{5-11} \cline{6-11} \cline{7-11} \cline{8-11} \cline{9-11} \cline{10-11} \cline{11-11} 
 & {\tiny{}rADL0} & {\tiny{}iALM} & {\tiny{}IPL} & {\tiny{}QP} & {\tiny{}RQP} & {\tiny{}rADL0} & {\tiny{}iALM} & {\tiny{}IPL} & {\tiny{}QP} & {\tiny{}RQP}\tabularnewline
\hline 
{\tiny{}100} & \textbf{\tiny{}388} & {\tiny{}66000} & {\tiny{}6863} & {\tiny{}37470} & {\tiny{}8293} & \textbf{\tiny{}4.4} & {\tiny{}323.3} & {\tiny{}68.7} & {\tiny{}344.6} & {\tiny{}85.6}\tabularnewline
{\tiny{}200} & \textbf{\tiny{}486} & {\tiny{}70551} & {\tiny{}6902} & {\tiny{}37696} & {\tiny{}1475} & \textbf{\tiny{}5.6} & {\tiny{}334.9} & {\tiny{}66.9} & {\tiny{}335.4} & {\tiny{}13.4}\tabularnewline
{\tiny{}400} & \textbf{\tiny{}674} & {\tiny{}72760} & {\tiny{}6902} & {\tiny{}37972} & {\tiny{}1562} & \textbf{\tiny{}7.6} & {\tiny{}347.5} & {\tiny{}67.9} & {\tiny{}339.0} & {\tiny{}14.2}\tabularnewline
{\tiny{}1600} & \textbf{\tiny{}1090} & {\tiny{}74200} & {\tiny{}6921} & {\tiny{}38203} & {\tiny{}1309} & {\tiny{}12.6} & {\tiny{}361.4} & {\tiny{}68.9} & {\tiny{}346.3} & \textbf{\tiny{}12.1}\tabularnewline
{\tiny{}3200} & {\tiny{}1400} & {\tiny{}74568} & {\tiny{}6921} & {\tiny{}38243} & \textbf{\tiny{}1327} & {\tiny{}16.0} & {\tiny{}369.8} & {\tiny{}74.1} & {\tiny{}352.3} & \textbf{\tiny{}12.1}\tabularnewline
\end{tabular}}
\par\end{centering}
\caption{Innermost iteration counts and runtimes for LCQM problems.\label{tab:qmp_iter_run}}
\end{table}

\begin{table}[tbh]
\begin{centering}
\makebox[\textwidth][c]{%
\begin{tabular}{>{\centering}m{0.7cm}|>{\centering}p{0.6cm}>{\centering}p{0.6cm}>{\centering}p{0.6cm}>{\centering}p{0.6cm}>{\centering}p{0.6cm}|>{\centering}p{0.6cm}>{\centering}p{0.6cm}>{\centering}p{0.6cm}>{\centering}p{0.6cm}>{\centering}p{0.6cm}}
{\scriptsize{}$M$} & \multicolumn{5}{c|}{{\scriptsize{}$c_{\max}$}} & \multicolumn{5}{c}{{\scriptsize{}$c_{{\rm wavg}}/c_{\max}$}}\tabularnewline
\cline{2-11} \cline{3-11} \cline{4-11} \cline{5-11} \cline{6-11} \cline{7-11} \cline{8-11} \cline{9-11} \cline{10-11} \cline{11-11} 
 & {\tiny{}rADL0} & {\tiny{}iALM} & {\tiny{}IPL} & {\tiny{}QP} & {\tiny{}RQP} & {\tiny{}rADL0} & {\tiny{}iALM} & {\tiny{}IPL} & {\tiny{}QP} & {\tiny{}RQP}\tabularnewline
\hline 
{\tiny{}100} & \textbf{\tiny{}4E+1} & {\tiny{}2E+3} & {\tiny{}6E+2} & {\tiny{}1E+3} & {\tiny{}1E+3} & {\tiny{}0.27} & {\tiny{}0.08} & {\tiny{}0.96} & {\tiny{}0.30} & {\tiny{}0.01}\tabularnewline
{\tiny{}200} & \textbf{\tiny{}8E+1} & {\tiny{}3E+3} & {\tiny{}1E+3} & {\tiny{}3E+3} & {\tiny{}3E+3} & {\tiny{}0.29} & {\tiny{}0.08} & {\tiny{}0.97} & {\tiny{}0.30} & {\tiny{}0.08}\tabularnewline
{\tiny{}400} & \textbf{\tiny{}2E+2} & {\tiny{}6E+3} & {\tiny{}3E+3} & {\tiny{}5E+3} & {\tiny{}5E+3} & {\tiny{}0.33} & {\tiny{}0.08} & {\tiny{}0.97} & {\tiny{}0.31} & {\tiny{}0.11}\tabularnewline
{\tiny{}1600} & \textbf{\tiny{}6E+2} & {\tiny{}2E+4} & {\tiny{}1E+4} & {\tiny{}2E+4} & {\tiny{}2E+4} & {\tiny{}0.39} & {\tiny{}0.08} & {\tiny{}0.97} & {\tiny{}0.31} & {\tiny{}0.12}\tabularnewline
{\tiny{}3200} & \textbf{\tiny{}1E+3} & {\tiny{}5E+4} & {\tiny{}2E+4} & {\tiny{}4E+4} & {\tiny{}4E+4} & {\tiny{}0.39} & {\tiny{}0.08} & {\tiny{}0.97} & {\tiny{}0.31} & {\tiny{}0.13}\tabularnewline
\end{tabular}}
\par\end{centering}
\caption{Penalty parameter statistics for LCQM problems.\label{tab:qmp_penalty}}
\end{table}

\subsection{Comments about Numerical Experiments}

Algorithm rADL0 is generally the most efficient in terms of total inner (or ACG) iterations, runtime, and final penalty parameter used. Moreover, the experiments in Subsection~\ref{subsec:qvp} demonstrate that the adaptivity of $m$ (or equivalently $\lam$) substantially improves AIDAL in terms of both inner (or ACG) iteration count and runtime. Finally, while the penalty ratio $c_{\rm wavg} / c_{\max}$ is generally the lowest for iALM, the performance for iALM in terms of the number of innermost iterations and runtime is generally the worst among the tested methods.

\section{Concluding Remarks}

\label{sec:concl_remarks}

Similar to the analyses in \cite{ImprovedShrinkingALM20,Lin2019},
the analysis of the AIDAL method strongly makes use of assumption
(A3) and the assumption that $D_{h}<\infty$ to obtain its competitive
${\cal O}(\varepsilon^{-5/2}\log\varepsilon^{-1})$ iteration complexity
when $\varepsilon=\rho=\eta$.
However, we conjecture that these two assumptions may be removed
using the more complicated analysis in \cite{RJWIPAAL2020} to obtain
a slightly worse ${\cal O}(\varepsilon^{-3}\log\varepsilon^{-1})$
iteration complexity (like in \cite{RJWIPAAL2020}). 

Like the adaptive prox-stepsize AIDAL in Appendix~\ref{app:adapt_aidal}, another
possible extension of AIDAL is one in which $\lam$, $\chi$, and $\theta$ are simultaneously chosen in an adaptive manner.
Moreover, it would be interesting to develop such an adaptive AIDAL and show that
it has the 
same iteration complexity bound
as the nonadaptive AIDAL in Algorithm~\ref{alg:aidal}.

\appendix

\section{Key Technical Bounds}

\label{app:mult_bd}

The appendix presents a key technical bound that is used in the analysis of AIDAL. 

\begin{lem}
\label{lem:tech_ADeltaZ_DeltaP}For every $(\tau,\theta)\in[0,1]^{2}$
satisfying $\tau\leq\theta^{2}$ and every $a,b\in\rn$, we have that
\begin{equation}
\|a-(1-\theta)b\|^{2}-\tau\|a\|^{2}\geq\left[\frac{(1-\tau)-(1-\theta)^{2}}{2}\right]\left(\|a\|^{2}-\|b\|^{2}\right).\label{eq:inexact_tech_DeltaP1}
\end{equation}
\end{lem}

\begin{proof}
Let $a,b\in\rn$ be fixed and define 
\begin{equation}
z=\left[\begin{array}{c}
\|a\|\\
\|b\|
\end{array}\right],\quad M=\left[\begin{array}{cc}
(1-\tau)+(1-\theta)^{2} & -2(1-\theta)\\
-2(1-\theta) & (1-\tau)+(1-\theta)^{2}
\end{array}\right].\label{eq:aux_tech_DeltaP1_def}
\end{equation}
Moreover, using our assumption of $\tau\leq\theta^{2}\le1$, observe
that 
\begin{align*}
\det M & =\left[(1-\tau)+(1-\theta)^{2}-2(1-\theta)\right]\left[(1-\tau)+(1-\theta)^{2}+2(1-\theta)\right]\\
 & =\left[\theta^{2}-\tau\right]\left[(1-\tau)+(1-\theta)^{2}+2(1-\theta)\right]\geq0,
\end{align*}
and hence, by Sylvester's criterion, it follows that $M\succeq0$.
Combining this fact with the Cauchy-Schwarz inequality and \eqref{eq:aux_tech_DeltaP1_def},
we thus have that 
\begin{align*}
 & \|a-(1-\theta)b\|^{2}-\tau\|a\|^{2}\geq(1-\tau)\|a\|^{2}-2(1-\theta)\|a\|\cdot\|b\|+(1-\theta)^{2}\|b\|^{2}\\
 & =\frac{1}{2}z^{T}Mz+\left[\frac{(1-\tau)-(1-\theta)^{2}}{2}\right]\left(\|a\|^{2}-\|b\|^{2}\right)\geq\left[\frac{(1-\tau)-(1-\theta)^{2}}{2}\right]\left(\|a\|^{2}-\|b\|^{2}\right). \qedhere
\end{align*}
\end{proof}

\section{Statement and Analysis of the ACG Algorithm}

\label{app:ACG}

Recall from Section~\ref{sec:intro} that our interest is in solving
\eqref{eq:main_prb} by inexactly solving NCO subproblems of the form
in \eqref{eq:primal_update}. This subsection presents an ACG algorithm
for inexactly solving latter type of problem and it considers the
more general class of NCO problems 
\begin{equation}
\min_{u\in\rn}\left\{ \psi(u) := \psi_{s}(u)+\psi_{n}(u)\right\} ,\label{eq:gen_nco}
\end{equation}
where the functions $\psi_{s}$ and $\psi_{n}$ are assumed to satisfy
the following assumptions: 
\begin{itemize}
\item[(B1)] $\psi_{n}:\rn\mapsto(-\infty,\infty]$ is a proper closed convex function.
\item[(B2)] $\psi_{s}$ is $\mu$-strongly convex and continuously differentiable on $\rn$ and satisfies
\begin{equation}
\label{eq:acg_lipsch_curv}
\|\nabla\psi_s(z) - \nabla\psi_s(z')\| \leq L\|z-z'\|
\end{equation}
for every $z',z\in\rn$ and some $L > 0$ and $\mu \in (0, L]$. 
\end{itemize}
Clearly, problem \eqref{eq:primal_update} is a special case of \eqref{eq:gen_nco},
and hence, any result that is stated in the context of \eqref{eq:gen_nco}
also applies to \eqref{eq:primal_update}. It is also well-known that assumption (B2) implies
\begin{equation}
\label{eq:acg_curv}
\frac{\mu}{2} \|z'-z\|^{2} \leq \psi_{s}(z')-\ell_{\psi_{s}}(z';z)\leq \frac{L}{2}\|z'-z\|^{2},
\end{equation}
for every $z,z'\in\rn$.

The pseudocode for the ACG algorithm is stated in Algorithm~\ref{alg:acg} which,
for a given a pair $({\sigma},x_{0})\in\r_{++}\times\dom\psi_{n}$,
inexactly solves \eqref{eq:gen_nco} by obtaining a pair $(z,v)$
satisfying 
\begin{equation}
v\in \nabla \psi_{s}(z)+ \partial\psi_{n}(z),\quad\|v\|\leq{\sigma}\|z-x_{0}\|.\label{eq:approx_acg_soln}
\end{equation}
Note that if ACG algorithm obtains the aforementioned triple with
${\sigma}=0$ then the first component of the triple is, in fact, a
global solution of \eqref{eq:gen_nco}.
 Indeed, if ${\sigma}=0$ then the above
inequality implies that $v=0$,
and the above inclusion reduces to
$0\in\partial(\psi_{s}+\psi_{n})(z)$,
which in view of \eqref{def:epsSubdiff} clearly implies that
$z$ is a global solution of
\eqref{eq:gen_nco}.

\begin{algorithm}[!htb] 
\caption{Accelerated Composite Gradient (ACG) Algorithm}
\label{alg:acg}

\Require{$({\sigma},x_{0})\in \R_{++}\times\dom\psi_{n}$.}

\Output{a pair $(z,v)\in\dom\psi_{n}\times\r^{n}$
satisfying \eqref{eq:approx_acg_soln}.}

\BlankLine
\Fn{\ACG{$\{\psi_{s},\psi_{n}\},\{L,\mu\},{\sigma},x_0$}}{

\texttt{\textcolor{blue}{STEP 0}}\textcolor{blue}{{} (initialization)}:\;

Set $y_{0}\gets x_{0}$, $A_{0}\gets0$.\;

\For{$j \gets 0,1,...$}{

\texttt{\textcolor{blue}{STEP 1}}\textcolor{blue}{{} (main iterates)}:

$\text{\textbf{find} the positive scalar } a_{j} \text{ satisfying } a_{j}^2 = \frac{(1+\mu A_{j})(a_{j} + A_{j})}{L}$

$A_{j+1}\gets A_{j}+ a_{j}$\;

$\tilde{x}_{j}\gets\frac{A_{j}}{A_{j+1}}x_{j}+\frac{A_{j+1}-A_{j}}{A_{j+1}}y_{j}$\;

$x_{j+1}\gets\argmin_{y\in\rn}\left\{ \ell_{\psi_s}(y;\tilde{x}_j)+\psi_{n}(y)+\frac{L+\mu}{2}\|y-\tilde{x}_j\|^{2}\right\} $\;

$y_{j+1} \gets y_{j} + \frac{a_{j}}{1+\mu A_{j+1}} [L(x_{j+1}-\tilde{x}_{j}) + \mu(x_{j+1} - y_{j})]$

\texttt{\textcolor{blue}{STEP 2}}\textcolor{blue}{{} (termination check)}:

$u_{j+1}\gets \nabla \psi_s(x_{j+1}) - \nabla \psi_s(\tilde{x}_j) + (L+\mu) (\tilde{x}_j - x_{j+1})$

\If{$\|u_{j+1}\|\leq{\sigma}\|x_{j+1}-x_{0}\|$}{

$(z,v) \gets (x_{j+1}, u_{j+1})$

\Return{$(z, v)$}

}

} 

} 

\end{algorithm}

{\color{purple}
We now devote the remainder of the section to proving the following properties about the ACG algorithm. Variations of the arguments that follow can also be found in \cite{kong2022complexity,sujanani20222}.}
\begin{prop}
\label{prop:acg_props}The following properties hold about the ACG
algorithm: 
\begin{itemize}
\item[(a)] for every $j\geq0$, it holds that 
\[
u_{j+1}\in \nabla \psi_{s}(x_{j+1}) + \partial \psi_{n}(x_{j+1}) = \pt(\psi_s + \psi_n)(x_{j+1});
\]
\item[(b)] it stops in a number of iterations bounded above by 
\begin{equation}
\color{purple}
\label{eq:acg_iter_bd}
\left\lceil 1 + 2\sqrt{\frac{L}{\mu}}\log_{1}^{+} \left\{ \frac{4L(L+\mu)^2}{\mu\sigma^2} \right\} \right\rceil,
\end{equation}
and its output $(z,v)$ satisfies \eqref{eq:approx_acg_soln}.
\end{itemize}
\end{prop}

{\color{purple}

We first present some technical properties about the generated iterates of Algorithm~\ref{alg:acg}.

\begin{lem} 
\label{lem:basic_acg_props}
Define the quantities
\begin{align}
\tau_{j} & := 1 + \mu A_{j}, \label{eq:tau_j} \\
\tilde{q}_{j+1}(\cdot) & := \ell_{\psi_s}(\cdot;\tilde{x}_j) + \psi_n(\cdot) + \frac{\mu}{2}\|\cdot-\tilde{x}_j\|^2 \label{eq:tilde_q_j}\\
q_{j+1}(\cdot) & := \tilde{q}_j(x_{j+1}) + L \langle \tilde{x}_j - x_{j+1}, \cdot - x_{j+1} \rangle + \frac{\mu}{2}\|\cdot - x_{j+1}\|^2, \label{eq:q_j}
\end{align}
for every $j\geq 0$. Then, for every $j\geq 1$, the following statements hold:
\begin{itemize}
    \item[(a)] $A_{j+1} \geq \left[1 + \sqrt{\mu}/(2\sqrt{L})\right]^{2j} / L$;
    \item[(b)] $x_{j+1} = \argmin_{x} \{{q}_{j+1}(x) + L\|x-\tilde{x}_j\|^2/2\}$;
    \item[(c)] $y_{j+1} = \argmin_y\{a_j q_{j+1}(y) + \tau_{j} \|y-y_j\|^2/2\}$;
    \item[(d)] $q_{j+1}(\cdot) \leq \psi(\cdot)$.
\end{itemize}
\end{lem}

\begin{proof}
(a) See, for example, \cite[Lemma~4]{MontSvaiter_fista}.

(b) Since $\nabla q_{j+1}(x_{j+1}) = L(\tilde{x}_j - x_{j+1})$, it follows that $x_{j+1}$ satisfies the optimality condition of the given minimization problem. Hence, the desired identity follows.

(c) It follows from the definition of ${q}_{j+1}(\cdot)$ and the update rule of $y_{j+1}$ that $a_j \nabla q_{j+1}(y_{j+1}) = \tau_{j+1} (y_{j+1}-y_j)$. The conclusion now follows from the optimality condition for the desired identity.

(d) In view of \eqref{eq:acg_curv} and the definition of $\tilde{q}_{j+1}$, we first have that $\tilde{q}_{j+1}(\cdot) \leq \psi(\cdot)$. On the other hand, it follows from the optimality condition of $\tilde{x}_{j+1}$ in Algorithm~\ref{alg:acg}, the convexity of $\psi_n$, and the definition of $q_j(\cdot)$ that $L(\tilde{x}_j - x_{j+1}) \in \partial \tilde{q}_{j+1}(x_{j+1})$. Furthermore, since $\tilde{q}_{j+1}$ is $\mu$-strongly convex, we also have $L(\tilde{x}_j - x_{j+1}) \in \partial (\tilde{q}_{j+1} - \mu\|\cdot-x_{j+1}\|^2/2)(x_{j+1})$. Combining all these facts with the definition of the subdifferential, we thus conclude that 
\[
\psi(\cdot) \geq \tilde{q}_{j+1}(\cdot) \geq 
\tilde{q}_{j+1}(x_{j+1}) + L\langle \tilde{x}_j - x_{j+1}, \cdot - x_{j+1} \rangle + \frac{\mu}{2}\|\cdot-x_{j+1}\|^2 = q_{j+1}(\cdot).
\qedhere
\]
\end{proof}

The next result establishes an important technical bound.

\begin{lem}
For every $j\geq 0$ and $y\in \rn$, it holds that 
\begin{align}
\label{eq:main_tech_acg_ineq}
\begin{gathered}
A_{j} q_{j+1}(x_j) + a_j q_{j+1}(y) + \frac{\tau_{j}}{2}\|y_j - y\|^2 - \frac{\tau_{j+1}}{2}\|y_{j+1}-y\|^2 \\
\quad \geq A_{j+1} \left[\psi(x_{j+1}) + \frac{\mu}{2}\|x_{j+1} - \tilde{x}_j\|^2 \right],
\end{gathered}
\end{align}
where $\tau_j$ and $q_j(\cdot)$ are as in \eqref{eq:tau_j} and \eqref{eq:q_j}, respectively.
\end{lem}

\begin{proof}
Using the update rule for $A_{j+1}$ we first note that $\tau_{j+1} = \tau_j + \mu a_j$.
Combining this fact, the optimality condition in Lemma~\ref{lem:basic_acg_props}(c) and the fact that $a_j q_{j+1}(\cdot) + \tau_{j}\|\cdot-y_j\|^2/2$ is $\tau_{j+1}$-strongly convex, we then have that
\begin{equation}
a_j q_{j+1}(y)+\frac{\tau_{j}}{2}\|y-y_j\|^2 - \frac{\tau_{j+1}}{2}\|y-y_{j+1}\|^2 \geq a_j q_{j+1}(y_{j+1})+\frac{\tau_{j}}{2}\|y_{j+1}-y_j\|^2
\label{eq:tech_acg_ineq1}
\end{equation}
for every $y\in\rn$. On the other hand, using the convexity of $q_{j+1}(\cdot)$, the second bound in \eqref{eq:acg_curv}, Lemma~\ref{lem:basic_acg_props}(b), and the quadratic subproblem associated with $a_j$, we have 
\begin{align}
 & A_{j}q_{j+1}(x_{j})+a_{j}q_{j+1}(y_{j+1})+\frac{\tau_{j}}{2}\|y_{j+1}-y_{j}\|^{2} \nonumber\\
 & \geq A_{j+1}q_{j+1}\left(\frac{A_{j}x_{j}+a_{j}y_{j+1}}{A_{j+1}}\right)+\frac{\tau_{j}A_{j+1}^{2}}{2a_{j}^{2}}\left\Vert \frac{A_{j}x_{j}+a_{j}y_{j+1}}{A_{j+1}}-\frac{A_{j}x_{j}+a_{j}y_{j}}{A_{j+1}}\right\Vert ^{2} \nonumber\\
 & \geq A_{j+1}\min_{x\in\mathbb{R}^{n}}\left\{ q_{j+1}(x)+\frac{\tau_{j}A_{j+1}^{2}}{2a_{j}^{2}}\|x-\tilde{x}_{j}\|^{2}\right\} =A_{j+1}\min_{x\in\mathbb{R}^{n}}\left\{ q_{j+1}(x)+\frac{L}{2}\|x-\tilde{x}_{j}\|^{2}\right\} \nonumber\\
 & =A_{j+1}\left[q_{j+1}(x_{j+1})+\frac{L}{2}\|x_{j+1}-\tilde{x}_{j}\|^{2}\right]\geq
A_{j+1} \left[\psi(x_{j+1}) +  \frac{\mu}{2}\|x_{j+1}-\tilde{x}_{j}\|^{2} \right]. \label{eq:tech_acg_ineq2}
\end{align}
The conclusion follows from combining \eqref{eq:tech_acg_ineq1} and \eqref{eq:tech_acg_ineq2}.
\end{proof}

We now derive a general telescopic bound on the quantity $\|x_{j+1} - \tilde{x}_j\|^2$.

\begin{lem}
For every $j\geq 0$ and $x\in\rn$,  it holds that
\begin{equation}
\label{eq:eta_bd}
    \frac{\mu A_{j+1}}{2}\|x_{j+1} - \tilde{x}_j\|^2 \leq \eta_j(x) - \eta_{j+1}(x),
\end{equation}
where the potential $\eta_i(\cdot)$ is given by
\begin{equation}
    \label{eq:acg_eta_bd}
    \eta_i(\cdot) := A_i [\psi(x_i) - \psi(\cdot)] + \frac{\tau_i}{2} \|\cdot - y_i\|^2 \quad \forall i\geq 0.
\end{equation}
\end{lem}

\begin{proof} Subtracting $A_{j+1} \psi(y)$ from \eqref{eq:main_tech_acg_ineq} and using Lemma~\ref{lem:basic_acg_props}(d), we have that
\begin{align*}
& \frac{A_{j+1}}{2}\|x_{j+1}-\tilde{x}_j\|^2 + A_{j+1}\left[\psi(x_{j+1}) - \psi(y)\right] \\ 
& \leq A_{j} q_{j+1}(x_j) + a_j q_{j+1}(y) - A_{j+1} \psi(y) + \frac{\tau_{j}}{2}\|y_j - y\|^2 - \frac{\tau_{j+1}}{2}\|y_{j+1}-y\|^2 \\
& \leq A_{j} \psi(x_j) + a_j \psi(y) - A_{j+1} \psi(y) + \frac{\tau_{j}}{2}\|y_j - y\|^2 - \frac{\tau_{j+1}}{2}\|y_{j+1}-y\|^2.
\end{align*}
The conclusion follows by re-arranging the above bound and using the update rule for $A_{j+1}$ and the definition of $\eta_i(\cdot)$.
\end{proof}

Specializing the above result, we establish a bound for the residuals $\{u_{j+1}\}_{j\geq 0}$ in terms of the prox residual $\|x_{j+1} - x_0\|^2$.

\begin{lem}
For every $j \geq 0$, it holds that 
\begin{equation}
\label{eq:acg_uj_bd}
\|u_{j+1}\|^2 \leq \frac{4(L+\mu)^2}{\mu A_{j+1}} \|x_{j+1}-x_0\|^2.
\end{equation}
\end{lem}

\begin{proof} 
Using assumption (B2), the definition of $u_{j+1}$, the bound $(a+b)^2 \leq 2a^2 + 2b^2$ for $a,b\in\r$, \eqref{eq:eta_bd} at $x=x_j$, and the fact that $(A_0,\tau_0)=(0,1)$, we have that
\begin{align*}
\frac{\mu A_{j+1}\|u_{j+1}\|^{2}}{2} & \leq\frac{\mu\sum_{i=0}^{j}A_{i+1}\|u_{i+1}\|^{2}}{2}\\
 & =\frac{\mu\sum_{i=0}^{j}A_{i+1}\|\nabla\psi_{s}(x_{i+1})-\nabla\psi_{s}(\tilde{x}_{i})+(L+\mu)(\tilde{x}_{i}-x_{i+1})\|^{2}}{2}\\
 & 
 \overset{\text{(B2)}}{\leq} \mu\sum_{i=0}^{j}A_{i+1}\left[\|\nabla\psi_{s}(x_{i+1})-\nabla\psi_{s}(\tilde{x}_{i})\|^{2}+(L+\mu)^2\|\tilde{x}_{i}-x_{i+1}\|^{2}\right]\\
 & \overset{\eqref{eq:eta_bd}}{\leq} 2\mu(L+\mu)^2\sum_{i=0}^{j}A_{i+1}\|\tilde{x}_{i}-x_{i+1}\|^{2}\leq4(L+\mu)^2\left[\eta_{0}(x_{j+1})-\eta_{k+1}(x_{j+1})\right]\\
 & \overset{(A_{0}, \tau_{0}) = (0,1)}{=}4(L+\mu)^2\left[\frac{1}{2}\|x_{0}-x_{j+1}\|^{2}-\frac{\tau_{j+1}}{2}\|x_{0}-x_{j+1}\|^{2}\right]\\
 & \leq2(L+\mu)^2\|x_{0}-x_{j+1}\|^{2}. \qedhere
\end{align*}
\end{proof}

We are now ready to prove Proposition~\ref{prop:acg_props}.

\begin{proof}[Proof of Proposition~\ref{prop:acg_props}]
(a) Using the optimality of $x_{j+1}$ the definition of $u_{j+1}$ in Algorithm~\ref{alg:acg}, we have that
\begin{align*}
0 & \in \nabla\psi_s({\tilde{x}_j}) + \pt 
\psi_n(x_{j+1}) + (L+\mu)(x_{j+1} - \tilde{x}_j)
= -u_{j+1} + \nabla\psi_s({x}_{j+1}) + \pt 
\psi_n(x_{j+1})\\
& = -u_{j+1} + \pt(\psi_s + \psi_n)(x_{j+1})
\end{align*}
where the last identity follows from the fact that $\psi_s$ and $\psi_n$ are convex (see (B1)--(B2)).

(b) Let $J$ denote the quantity in \eqref{eq:acg_iter_bd}. Using Lemma~\ref{lem:basic_acg_props}(a) and the bound $\log(1+t) \geq t/2$ for $t\in[0,1]$, it is straightforward to verify that 
$4(L+\mu)^2/(\mu A_{J+1}) \leq \sigma^2$.
It then follows from the previous bound and \eqref{eq:acg_uj_bd} that 
\[
\|u_{J+1}\|^2 \leq  \frac{4(L+\mu)^2}{\mu A_{J+1}} \|x_{J+1}-x_0\|^2 \leq \sigma^2 \|x_{J+1}-x_0\|^2.
\]
Consequently, it follows from the above bound, part (a), and the termination condition of Algorithm~\ref{alg:acg} that the ACG algorithm stops in a number of iterations bounded above by $J$.
\end{proof}}

\section{Necessary Optimality Conditions}
\label{app:local_necessary}

This appendix shows that if $\hat{z}$ local minimum of \eqref{eq:main_prb} then condition \eqref{eq:stationary_soln} holds. Throughout this appendix, we denote
\[
\psi'(x;d) = \lim_{t\downarrow 0} \frac{\psi(x+td)-\psi(x)}{t}
\]
as the directional derivative of a function $\psi$ at $x$ in the direction $d$.

The first useful result presents a relationship between directional derivatives of composite functions and the usual first-order necessary conditions.

\begin{lem}
\label{lem:tech_ddir_min}
Let $g:\rn\mapsto(-\infty,\infty]$ be a proper convex function, and let $f$ be a differentiable function on $\dom g$. Then, for every $x\in\dom g$, the following statements hold:
\begin{itemize}
    \item[(a)] $\inf_{\|d\|\leq 1} (f+g)'(x;d) = -\inf_{u\in\rn} \{\|u\| : u \in\nabla f(x) +\partial g(x)\}$;
    \item[(b)] if $x$ is a local minimum of $f+h$ then {\color{purple}$0 \in \nabla f(x)+\partial h(x)$}.
\end{itemize}
\end{lem}

\begin{proof}
(a) See \cite[Lemma 15]{Kong2019} with $({\cal X}, h)=(\rn, g)$.

(b) This follows immediately from (a) and the fact that $(f+h)'(x;d)\geq 0$ for every $d\in\rn$.
\end{proof}

We now establish the aforementioned necessary condition.

\begin{prop} 
Let $(f,h,A,b)$ be as in (A1)-(A4). If $\hat{z}$ is a local minimum of \eqref{eq:main_prb}, then there exists a multiplier $\hat{p}$ such that \eqref{eq:stationary_soln} holds.
\end{prop}

\begin{proof}

We first establish an important technical identity. Let $S=\{z \in \rn : Az=b\}$, let $\delta_S$ denote the indicator function of $S$, i.e., the function that takes value $0$ if its input is in $S$ and $+\infty$ otherwise, and let $\ri X$ denote the relative interior of a set $X$. Since assumptions (A3)--(A4) imply that $\ri{\cal H} \cap \ri{S} = \intr {\cal H} \cap {S} \neq \emptyset$, it follows from \cite[Theorem 23.8]{Rockafellar70} that for every $x\in {\cal H} \cap S$ we have
\begin{equation}
\partial(\delta_S + h)(x) = \partial \delta_S(x) + \partial h(x) = N_S(x) + \partial h(x) = \{\xi + A^*p : \xi \in \partial h(x)\}. \label{eq:tech_local_subdiff}
\end{equation}
The conclusion follows from the above identity and Lemma~\ref{lem:tech_ddir_min}(b) with $g=h+\delta_{S}$.
\end{proof}

\section{Adaptive AIDAL}

\label{app:adapt_aidal}

This appendix presents an adaptive version of AIDAL where we choose the prox stepsize adaptively. 

Before presenting the algorithm, we first motivate its construction under the assumption that the reader is familiar with the notation and results of Section~\ref{sec:cvg_analysis}. 
To begin, the careful reader may notice that the special choice of $\lam=1/(2m)$ in AIDAL (Algorithm~\ref{alg:aidal}) is only needed to ensure that the function $\lam {\cal L}_c^\theta(\cdot;p) + \|\cdot\|^2$ is strongly convex with respect to the norm $\|x\|_Q
 = \langle x, [(1-\lam m)I + c\lam A^*A]x \rangle$ for every $c>0$ and $p\in A(\rn)$. Moreover, this global property is only needed to show that:
\begin{itemize}
    \item[(i)] the $k^{\rm th}$ ACG call of AIDAL stops with a pair $(z_k, v_k)$ satisfying $\|v_k\|\leq \sigma \|z_k - z_{k-1}\|$;
    \item[(ii)] $\lam \|\hat{v}_i\| \apprle \Psi_{k-1}^\theta - \Psi_{k}^\theta$.
\end{itemize}
The other technical details of Section~\ref{sec:cvg_analysis}, such as the boundedness of $\Psi_i^\theta$, are straightforward to show as long as the prox stepsize is bounded.
As a consequence, a natural relaxation of AIDAL is to employ a line search at its $k^{\rm th}$ outer iteration for the largest $\lam$ within a bounded range satisfying conditions (i) and (ii) above.

In Algorithm~\ref{alg:adap_aidal}, we present one possible relaxation. Specifically, the $k^{\rm th}$ prox stepsize $\lam_k$ is chosen from a set of candidates in the range $(0, \lam_{k-1}]$.

\begin{algorithm}[!htb] 
\caption{Adaptive AIDAL Method}
\label{alg:adap_aidal}

\Require{Same as in Algorithm~\ref{alg:aidal} but with additional parameters $\gamma > 1$ and $\lam_0 > 0$.}

\Output{Same as in Algorithm~\ref{alg:aidal}.}

\BlankLine
\Fn{\AdapAIDAL{$M$,$\{\sigma,\chi,\theta,\lam_0\},\{c_{1},z_{0},p_{0}\},\{{\rho},{\eta}\}$, $\gamma$}}{

$\lam_{0} \gets \lam$

\For{$k \gets 1,2,...$}{
\textbf{find} the smallest nonnegative integer $\beta_k$ such that the ACG call in step~1 of Algorithm~\ref{alg:aidal} with $\lam = \gamma^{-\beta_k} \lam_{k-1}$ stops with a pair $(z_k, v_k)$ satisfying 
\begin{equation}
\label{eq:adapt_aidal_cond}
\begin{cases}
\|v_{k}\|\leq\sigma\|z_{k}-z_{k-1}\| & \text{if }k\geq1, \text{\textbf{ and }}\\
\|v_{k}+z_{k-1}-z_{k}\|^{2}\leq 9 \lambda(\Psi_{k-1}^{\theta}-\Psi_{k}^{\theta}) & \text{if }k\geq2,
\end{cases}
\end{equation}
where $\Psi_k^\theta$ is given in \eqref{eq:Psi_def}

\textbf{set} $\lam_k \gets \gamma^{-\beta_k} \lam_{k-1}$

\textbf{execute} steps 1--4 of Algorithm~\ref{alg:aidal} with $\lam = \lam_k$

} 

} 

\end{algorithm}

We now make a few remarks about Algorithm~\ref{alg:adap_aidal}. First, the candidate search space for the $k^{\rm th}$ prox stepsize forms a geometrically decreasing sequence and $\lam_k \leq \lam_{k-1}$. Second, the first condition of \eqref{eq:adapt_aidal_cond} corresponds to condition (i), while the second condition corresponds to condition (ii). Moreover, the second condition of \eqref{eq:adapt_aidal_cond} always holds when $\lam = 1/(2m)$ due to Lemma~\ref{lem:alm_aux2}, Lemma~\ref{lem:alm_aux3}, and the definition of $\hat{v}_i$ which imply (cf. the proof of Proposition~\ref{prop:vi_poten_bd}) that
\begin{align*}
\|v_{k}+z_{k-1}-z_{k}\|^{2} & = \lambda^{2}\|\hat{v}_{k}\|^{2}\leq  9\lambda(\Psi_{k-1}^{\theta}-\Psi_{k}^{\theta}).
\end{align*}
Third, in view of the previous remark, since conditions (i) and (ii) are always satisfied whenever $\lam \leq 1/(2m)$, we also have that $\lam_k \in [1/(2\gamma m), \lam_0]$ and, hence, the sequence $\{\lam_k\}_{k\geq 1}$ is bounded.

Notice that it is not immediately clear how one obtains $\beta_k$ at the $k^{\rm th}$ outer iteration. One possible approach is to apply an adaptive ACG variant to the stepsize sequence $\{\lam_{k-1}\beta^{-j}\}_{j\geq0}$ in which the variant has a mechanism to determine if at least one of the conditions in \eqref{eq:adapt_aidal_cond} is reachable. 
This is so that if none of the conditions in \eqref{eq:adapt_aidal_cond} are reachable for some candidate $\lam$, then the variant can be called again with a smaller stepsize. One example is the adaptive ACG variant in \cite{kong2022complexity}, which contains a mechanism for determining the reachability of the first condition in \eqref{eq:adapt_aidal_cond} and can even adaptively choose its other curvature parameters, such as $L$ in Algorithm~\ref{alg:acg}. Note that if the ACG has already been called with the $\beta_k$ satisfying \eqref{eq:adapt_aidal_cond} during the $\beta_k$ line search, then it does not need to be called again when executing the steps of Algorithm~\ref{alg:aidal}.

Before closing this section, we briefly discuss the convergence and iteration complexity of the method. 
Convergence of the method is straightforward to establish using the same techniques of Section~\ref{sec:cvg_analysis} and the fact that $\lam_k$ is bounded (see the remarks above). 
On the other hand, it can be shown that the iteration complexity of the method is on the same order of complexity as in Theorem~\ref{thm:aidal_compl}. 
Without going through the cumbersome technical details, we assert that this follows from the boundedness of the stepsizes $\lam_k$, the fact that the search for the next stepsize is done geometrically, and arguments similar to other adaptive augmented Lagrangian/penalty methods such as the one in \cite{WJRVarLam2018}.

\subsection*{Data Availability Statement}

The data and code generated, used, and/or analyzed during the current
study are publicly available in the \texttt{NC-OPT} GitHub repository\footnote{See \href{https://github.com/wwkong/nc_opt}{https://github.com/wwkong/nc\_opt}.}
under the directory \texttt{./tests/papers/aidal/}.

\subsection*{Ethics Statement}

The authors declare that they have no conflict of interest.

{\small{}\bibliographystyle{plain}
\bibliography{Proxacc_ref}

\begin{thebibliography}{10}

\bibitem{Aybatpenalty}
N.~S. Aybat and G.~Iyengar.
\newblock A first-order smoothed penalty method for compressed sensing.
\newblock {\em SIAM J. Optim.}, 21(1):287--313, 2011.

\bibitem{AybatAugLag}
N.~S. Aybat and G.~Iyengar.
\newblock A first-order augmented {Lagrangian} method for compressed sensing.
\newblock {\em SIAM J. Optim.}, 22(2):429--459, 2012.

\bibitem{boob2022stochastic}
D.~Boob, Q.~Deng, and G.~Lan.
\newblock Stochastic first-order methods for convex and nonconvex functional
  constrained optimization.
\newblock {\em Math. Program.}, pages 1--65, 2022.

\bibitem{Goncalves2019}
M.~L.~N. Goncalves, J.~G. Melo, and R.~D.~C. Monteiro.
\newblock Convergence rate bounds for a proximal {ADMM} with over-relaxation
  stepsize parameter for solving nonconvex linearly constrained problems.
\newblock {\em Pac. J. Optim.}, 15(3):379--398, 2019.

\bibitem{NIPS2014_5615}
Q.~Gu, Z.~Wang, and H.~Liu.
\newblock Sparse {PCA} with oracle property.
\newblock In Z.~Ghahramani, M.~Welling, C.~Cortes, N.~D. Lawrence, and K.~Q.
  Weinberger, editors, {\em Adv. Neural Inf. Process. Syst. 27}, pages
  1529--1537. Curran Associates, Inc., 2014.

\bibitem{HongPertAugLag}
D.~Hajinezhad and M.~Hong.
\newblock Perturbed proximal primal-dual algorithm for nonconvex nonsmooth
  optimization.
\newblock {\em Math. Program.}, 176:207--245, 2019.

\bibitem{SZhang-Pen-admm}
B.~Jiang, T.~Lin, S.~Ma, and S.~Zhang.
\newblock Structured nonconvex and nonsmooth optimization algorithms and
  iteration complexity analysis.
\newblock {\em Comput. Optim. Appl.}, 72(3):115--157, 2019.

\bibitem{kong2021thesis}
W.~{Kong}.
\newblock {Accelerated Inexact First-Order Methods for Solving Nonconvex
  Composite Optimization Problems}.
\newblock {\em Available on arXiv:2104.09685}, April 2021.

\bibitem{kong2022complexity}
W.~Kong.
\newblock Complexity-optimal and curvature-free first-order methods for finding
  stationary points of composite optimization problems.
\newblock {\em arXiv preprint arXiv:2205.13055}, 2022.

\bibitem{WJRproxmet1}
W.~Kong, J.~G. Melo, and R.~D.~C. Monteiro.
\newblock Complexity of a quadratic penalty accelerated inexact proximal point
  method for solving linearly constrained nonconvex composite programs.
\newblock {\em SIAM J. Optim.}, 29(4):2566--2593, 2019.

\bibitem{WJRVarLam2018}
W.~Kong, J.~G. Melo, and R.~D.~C. Monteiro.
\newblock An efficient adaptive accelerated inexact proximal point method for
  solving linearly constrained nonconvex composite problems.
\newblock {\em Comput. Optim. Appl.}, 76(2):305--346, 2020.

\bibitem{kong2020iteration}
W.~Kong, J.~G. Melo, and R.~D.~C. Monteiro.
\newblock {Iteration-complexity of a proximal augmented {L}agrangian method for
  solving nonconvex composite optimization problems with nonlinear convex
  constraints}.
\newblock {\em Available on arXiv:2008.07080}, 2020.

\bibitem{Melo2020}
W.~Kong, J.~G. Melo, and R.~D.~C. Monteiro.
\newblock Iteration complexity of an inner accelerated inexact proximal
  augmented {L}agrangian method based on the classical {L}agrangian function.
\newblock {\em SIAM Journal on Optimization}, 33(1):181--210, 2023.

\bibitem{Kong2019}
W.~Kong and R.~D.~C. Monteiro.
\newblock An accelerated inexact proximal point method for solving
  nonconvex-concave min-max problems.
\newblock {\em SIAM J. Optim.}, 31(4):2558--2585, 2021.

\bibitem{LanRen2013PenMet}
G.~Lan and R.~D.~C. Monteiro.
\newblock Iteration-complexity of first-order penalty methods for convex
  programming.
\newblock {\em Math. Program.}, 138(1):115--139, Apr 2013.

\bibitem{LanMonteiroAugLag}
G.~Lan and R.~D.~C. Monteiro.
\newblock Iteration-complexity of first-order augmented {L}agrangian methods
  for convex programming.
\newblock {\em Math. Program.}, 155(1):511--547, Jan 2016.

\bibitem{ImprovedShrinkingALM20}
Z.~Li, P.-Y. Chen, S.~Liu, S.~Lu, and Y.~Xu.
\newblock Rate-improved inexact augmented {L}agrangian method for constrained
  nonconvex optimization.
\newblock {\em Int. Conf. Artif. Intell. Stat.}, pages 2170--2178, 2021.

\bibitem{Li2019}
Z.~Li and Y.~Xu.
\newblock Augmented {L}agrangian--based first-order methods for
  convex-constrained programs with weakly convex objective.
\newblock {\em INFORMS Journal on Optimization}, 3(4):373--397, 2021.

\bibitem{Lin2019}
Q.~Lin, R.~Ma, and Y.~Xu.
\newblock Inexact proximal-point penalty methods for constrained non-convex
  optimization.
\newblock {\em Available on arXiv:1908.11518}, 2019.

\bibitem{ShiqiaMaAugLag16}
Y.-F. Liu, X.~Liu, and S.~Ma.
\newblock On the nonergodic convergence rate of an inexact augmented
  {L}agrangian framework for composite convex programming.
\newblock {\em Mathematics of Operations Research}, 44(2):632--650, 2019.

\bibitem{zhaosongAugLag18}
Z.~{Lu} and Z.~{Zhou}.
\newblock {Iteration-complexity of first-order augmented Lagrangian methods for
  convex conic programming}.
\newblock {\em Available on arXiv:1803.09941}, 2018.

\bibitem{RJWIPAAL2020}
J.~G. Melo, R.~D.~C. Monteiro, and H.~Wang.
\newblock Iteration-complexity of an inexact proximal accelerated augmented
  {L}agrangian method for solving linearly constrained smooth nonconvex
  composite optimization problems.
\newblock {\em Available on arXiv:2006.08048}, 2020.

\bibitem{MontSvaiter_fista}
R.~D.~C. Monteiro, C.~Ortiz, and B.~F. Svaiter.
\newblock An adaptive accelerated first-order method for convex optimization.
\newblock {\em Comput. Optim. Appl.}, 64:31--73, 2016.

\bibitem{IterComplConicprog}
I.~Necoara, A.~Patrascu, and F.~Glineur.
\newblock Complexity of first-order inexact {L}agrangian and penalty methods
  for conic convex programming.
\newblock {\em Optim. Methods Softw.}, pages 1--31, 2017.

\bibitem{Patrascu2017}
A.~Patrascu, I.~Necoara, and Q.~Tran-Dinh.
\newblock Adaptive inexact fast augmented {L}agrangian methods for constrained
  convex optimization.
\newblock {\em Optim. Lett.}, 11(3):609--626, 2017.

\bibitem{Rockafellar70}
R.T. Rockafellar.
\newblock {\em Convex Analysis}.
\newblock Princeton University Press, Princeton, 1970.

\bibitem{inexactAugLag19}
M.~Sahin, A.~Eftekhari, A.~Alacaoglu, F.~Latorre, and V.~Cevher.
\newblock An inexact augmented {L}agrangian framework for nonconvex
  optimization with nonlinear constraints.
\newblock {\em Adv. Neural Inf. Process. Syst.}, 32, 2019.

\bibitem{sujanani20222}
A.~{Sujanani} and R.~D.~C. {Monteiro}.
\newblock {An adaptive superfast inexact proximal augmented Lagrangian method
  for smooth nonconvex composite optimization problems}.
\newblock {\em arXiv e-prints}, page arXiv:2207.11905, July 2022.

\bibitem{Xu2019}
Y.~Xu.
\newblock Iteration complexity of inexact augmented {L}agrangian methods for
  constrained convex programming.
\newblock {\em Math. Program.}, 2019.

\bibitem{ErrorBoundJzhang-ZQLuo2020}
J.~Zhang and Z.-Q. Luo.
\newblock A global dual error bound and its application to the analysis of
  linearly constrained nonconvex optimization.
\newblock {\em Available on arXiv:2006.16440}, 2020.

\bibitem{ADMMJzhang-ZQLuo2020}
J.~Zhang and Z.-Q. Luo.
\newblock A proximal alternating direction method of multiplier for linearly
  constrained nonconvex minimization.
\newblock {\em SIAM J. Optim.}, 30(3):2272--2302, 2020.

\end{thebibliography}
}{\small\par}
\end{document}